\documentclass[12pt,leqno]{amsart}
\usepackage{a4wide}
\usepackage{amssymb}
\usepackage{amsthm}
\usepackage{amsmath}
\usepackage{amscd}
\usepackage[mathscr]{eucal}
\usepackage[all]{xy}

\numberwithin{equation}{section}

\theoremstyle{plain}
\newtheorem{theorem}{Theorem}[section]
\newtheorem{corollary}[theorem]{Corollary}
\newtheorem{lemma}[theorem]{Lemma}
\newtheorem{proposition}[theorem]{Proposition}

\theoremstyle{definition}

\newtheorem{remark}[theorem]{Remark}
\newtheorem{example}[theorem]{Example}

\theoremstyle{remark}

\newcommand{\R}{\mathbb{R}}
\newcommand{\Q}{\mathbb{Q}}
\newcommand{\Z}{\mathbb{Z}}

\newcommand{\C}{\mathbb{C}}
\newcommand{\h}{\mathbb{H}}
\renewcommand{\H}{\mathbb{H}}

\newcommand{\G}{\Gamma}
\newcommand{\g}{\gamma}
\newcommand{\la}{\lambda}
\newcommand{\La}{\Lambda}

\newcommand{\back}{\backslash}

\newcommand{\calD}{\mathcal{D}}
\newcommand{\calE}{\mathcal{E}}
\newcommand{\calF}{\mathcal{F}}

\newcommand{\calL}{\mathcal{L}}

\newcommand{\calO}{\mathcal{O}}

\newcommand{\calZ}{\mathcal{Z}}

\newcommand{\frakc}{\mathfrak c}

\newcommand{\frake}{\mathfrak e}

\newcommand{\eps}{\varepsilon}

\newcommand{\vol}{\operatorname{vol}}

\newcommand{\sgn}{\operatorname{sgn}}

\newcommand{\Span}{\operatorname{span}}

\newcommand{\Sl}{\operatorname{SL}}

\newcommand{\Symp}{\operatorname{Sp}}
\newcommand{\Mp}{\operatorname{Mp}}
\newcommand{\Orth}{\operatorname{O}}
\newcommand{\Uni}{\operatorname{U}}

\newcommand{\SO}{\operatorname{SO}}

\newcommand{\dv}{\operatorname{div}}
\newcommand{\Ei}{\operatorname{Ei}}

\begin{document}

\title{On Two Geometric Theta Lifts}
\author[Jan H.~Bruinier and Jens Funke]{Jan Hendrik Bruinier* and Jens Funke**}
\date{March 14, 2003}
\address{Mathematisches Institut, Universit\"at zu K\"oln, Weyertal 86--90, D-50931 K\"oln, Germany}
\email{bruinier@math.uni-koeln.de } 
\address{Fields Institute, 222 College Street, Toronto, Ontario, M5T 3J1, Canada}
\email{jfunke@fields.utoronto.ca}

\thanks{* Supported by a Heisenberg-Stipendium of the DFG}
\thanks{** J.E.~Marsden PostDoctoral Fellow 02/03 at the Fields Institute, Toronto}

\maketitle


\section{Introduction}
\label{sect1}

Borcherds \cite{Bo1} introduced a `singular' theta lift from modular forms of (typically) negative weight for $\Sl_2(\R)$ to the orthogonal group $\Orth(p,2)$. This gave rise to remarkable product expansions of automorphic forms for $\Orth(p,2)$, but also to the realization of generating series of certain `special' divisors in locally symmetric spaces attached to $\Orth(p,2)$ as holomorphic modular forms of positive weight \cite{Bo2}. The Borcherds lift was generalized by the first named author and studied in connection with the cohomology of the special divisors \cite{Br1,Br2}. 

On the other hand, the Weil representation and theta series have been utilized in a more classical way by several people (e.g.~\cite{Sh,Od,RS,TW,W}) to establish generating series of special cycles in orthogonal and unitary locally symmetric spaces of arbitrary signature as modular forms, in greatest generality by Kudla and Millson (see e.g.~\cite{KM90}).

The purpose of this paper is threefold: 
\begin{itemize}
\item[(1)]
For $\Orth(p,2)$, we derive an adjointness result between the Kudla-Millson lift and the Borcherds lift.

\item[(2)]
We also introduce a new Borcherds lift for $\Orth(p,q)$ and obtain a similar relationship to the Kudla-Millson lift.

\item[(3)]
As a geometric application, we show that the exterior derivative of the current induced by this generalized Borcherds lift is closely related to the delta current of a special cycle. 

\end{itemize}

\medskip

We sketch the main results of the paper in more detail.

Let $V$ be a quadratic space over $\Q$ of signature $(p,q)$ and write $D$ for  the  associated symmetric space. Let $L$ be an even lattice in $V$, $\G$ be a  subgroup of finite index of the group of units of $L$ and $X= \G \back D$ the associated locally symmetric space. The special cycles arise from subsymmetric spaces of codimension $q$ induced by embeddings of orthogonal groups of signature $(p-1,q)$ into $\Orth(V)$.

For simplicity of the exposition we assume in the introduction that $L$ is unimodular. (In the main body of the paper, we treat the general case of arbitrary level by using Borcherds' vector valued modular form setting. Moreover, for $q$ odd, this is essential in order to obtain a nonzero theory).



We first introduce a new space of automorphic forms of (typically negative) weight $k$.  Its importance lies in the fact that on one hand we will  systematically treat this space as  the input space for both the Borcherds lift for $\Orth(p,2)$ and its extension to $\Orth(p,q)$, while on the other hand it satisfies a duality with the space of holomorphic cusp forms of (positive) weight $2-k$. This duality, interesting in its own right, will be crucial for us.

Namely, we let $H_k$  be the space of  \emph{weak Maass forms}, consisting of those smooth functions $f$ on the upper half plane $\h$, which transform with weight $k= 2-(p+q)/2$ under $\Sl_2(\Z)$, are annihilated by the Laplacian of weight $k$, and satisfy $f(\tau) = O(e^{Cv})$ as $\tau= u+ iv \to i\infty$ for some constant $C>0$. 
For $f \in H_k$, put $\xi_k(f) = R_{-k}(v^k\bar{f})$, where $R_{-k}$ is the standard raising operator for modular forms of weight $-k$. We prove that $\xi_k$ defines an antilinear map $\xi_k: H_k \to M^!_{2-k}$, to the space of 
meromorphic modular forms of weight $2-k$ which are holomorphic on $\H$ (see Proposition \ref{defmap}). 
It is easily checked that $M^!_{k}$ is the kernel of $\xi_k$. 
We let $H_k^+$ be the preimage of $S_{2-k}$, the space of cusp forms of weight $2-k$.

\begin{theorem}\label{IT1} 
The bilinear pairing between $S_{2-k}$ and $H^+_k$, defined by 
\[
\{g,f\} = (g, \xi_k(f))_{2-k},
\]
for $g \in S_{2-k}$ and $f \in H_k^+$, 
induces a non-degenerate pairing of $S_{2-k}$ and $H^+_k/M_k^!$.
Here $(\,,\,)_{2-k}$ is the Petersson scalar product for modular forms of weight $2-k$. 
\end{theorem}

The main point here is to show that the map $\xi_k$ is surjective. The proof is using methods from complex geometry and is ultimately based on Serre duality.
The pairing $\{g,f\}$ can be explicitly evaluated in terms of the Fourier coefficients of $g$ and the singular part of $f$ (Proposition \ref{proppair}).
Note that Borcherds \cite{Bo2} established a similar duality statement (in terms of formal power series). 

\smallskip

For signature $(p,2)$, when $D$ is Hermitian, we then introduce the Borcherds lift as a map on $H_k^+$; i.e, for any $f \in H_k^+$, its lift is given  by integrating $f \in H_k^+$ against $\Theta(\tau,z,\varphi_0)$, ($\tau=u+iv \in \h$, $z \in D$), the Siegel theta series attached to the standard Gaussian $\varphi_0$ on $V(\R)$ (see also \cite{Bo1}, \cite{Br1}):
\begin{equation}\label{IBlift0}
\Phi(z,f) = \int^{reg}_{\Sl_2(\Z) \back \h} f(\tau) \Theta(\tau,z,\varphi_0) \frac{du \, dv}{v^2}.
\end{equation}
The integral is usually divergent, and a suitable regularization was found by Harvey and Moore \cite{HM}. The regularization process leads to logarithmic singularities along certain special cycles $Z(f)$ in $D$ which in this case are divisors. Moreover, $\Lambda_B(z,f) := dd^c \Phi(z,f)$ actually extends to a closed smooth $(1,1)$-form on $X$. Here $d$ and $d^c$ are the standard (exterior) differential operators on $D$. Hence, we have a map
\begin{equation}\label{IBlift}
\La_B: H^+_k  \longrightarrow \mathcal{Z}^{2}(X), 
\end{equation}
to the space of closed $2$-forms on $X$. For $f \in M^!_k$ one has $\Lambda_B(f) = a^+(0) \Omega$, where $a^+(0)$ is the constant coefficient of $f$, and $\Omega$ is the suitably normalized K\"ahler form on $D$.

\smallskip

On the other hand, Kudla and Millson \cite{KM90} construct for general signature $(p,q)$ a theta series 
$\Theta(\tau,z,\varphi_{KM})$ associated to a certain Schwartz function $\varphi_{KM}$ on $V(\R)$, which takes values in $\mathcal{Z}^{q}(X)$. Then for $\eta \in 
\mathcal{Z}_c^{(p-1)q}(X)$, the compactly supported closed $(p-1)q$-forms on $X$, the Kudla-Millson lift is defined by 
\[
\Lambda_{KM}(\tau,\eta) = \int_X \eta \wedge \Theta(\tau,\varphi_{KM}).
\]
It turns out that $\Lambda_{KM}(\tau,\eta)$ is actually a holomorphic modular form of weight $2-k$, so that we have a map
\begin{equation}\label{IKMlift}
\Lambda_{KM}: \mathcal{Z}_c^{(p-1)q}(X) \longrightarrow M_{2-k},
\end{equation}
which also factors through cohomology. Moreover, the Fourier coefficients of $\Lambda_{KM}$ are given by periods of $\eta$ over the special cycles.


\begin{theorem}\label{IT2}
Assume $D$ be Hermitian, i.e., $q=2$, and 
let $f \in H^+_k$ with constant coefficient $a^+(0)$. We then have the following identity of closed $2$-forms on $X$:
\[
\Lambda_B(z,f) = (\Theta(\tau,z,\varphi_{KM}),\xi_k(f))_{2-k} + a^+(0) \Omega.
\]
Therefore the maps  $\Lambda_B$ 
and $\Lambda_{KM}$ 
are naturally adjoint via the standard pairing $(\,,\,)_X$ of 
$\mathcal{Z}_c^{2p-2}(X)$ with $\mathcal{Z}^2(X)$ and the paring $\{\,,\,\}$ of $M_{2-k}$ with $H^+_k$, that is,
\[
\big(\eta, \Lambda_{B}(f)\big)_X  =\left\{\Lambda_{KM}(\eta),f\right\} +   a^+(0) \big(\eta,\Omega\big)_X.
\]
Furthermore, this duality factors through cohomology, and $H^+_k/M_k^!$, respectively. (See also Theorem \ref{Main2}.)
\end{theorem}

This is based on the fundamental relationship between the two theta series involved:
\begin{theorem}\label{IT3}
Let $L_{2-k}$ be the lowering Maass operator of weight $2-k$ on $\h$. Then
\[
L_{2-k} \Theta(\tau,z,\varphi_{KM}) = -dd^c \Theta(\tau,z,\varphi_0).
\]
\end{theorem}
We show this by switching to the Fock model of the Weil representation. Then the idea for the proof of Theorem \ref{IT2} is given by the following formal (!) calculation:
\begin{align*}
\left(\Theta(\varphi_{KM}),\xi_k(f)\right)_{2-k}  
&\,\text{``$=$''}\, 
- \left(L_{2-k} \Theta(\varphi_{KM}), v^k \bar{f}\right)^{reg}_{-k} \\
&\,\text{``$=$''}\,  
 \left(dd^c \Theta(\varphi_0),v^k\bar{f}\right)^{reg}_{-k}   
 \,\text{``$=$''}\, dd^c \int^{reg}_{\Sl_2(\Z) \back \h} \Theta(\varphi_0) f \frac{dudv}{v^2}. 
\end{align*}
The first equality follows from the adjointness of the raising and lowering operator. However, the second scalar product no longer converges and needs to be regularized. At this point, one also obtains an `error' term involving the K\"ahler form. The second equality follows from the key fact provided in Theorem \ref{IT3} and the last one by interchanging the integration and differentiation (which also needs careful consideration). An additional difficulty occurs by the fact that in the intermediate steps one also has to deal with singularities along the special divisors. 

As an application, we are able to recover several of the geometric properties of the Borcherds lift in \cite{Br1} with often simpler proofs.

\smallskip

For general signature $(p,q)$, we note that by the 
work of Kudla and Millson \cite{KM90}, there exists a theta function $\Theta(\tau,z,\psi)$ with values in the $(q-1)$-forms of $X$ such that
\begin{equation*}
L_{2-k} \varphi_{KM}(\tau,z,\varphi_{KM}) = d \Theta(\tau,z,\psi).
\end{equation*}
 This is the exact analogon for the crucial identity given in Theorem \ref{IT3} in the Hermitian case. So, for $f \in H_k^+$, we define the Borcherds lift 
$\Phi(z,f,\psi)$ by
replacing $\Theta(\tau,z,\varphi_0)$  by $\Theta(\tau,z,\psi)$ 
as the kernel function in (\ref{IBlift0}). Note that this is now not scalar-valued but a $(q-1)$-form with singularities of higher order along the special cycle $Z(f)$ (see Proposition \ref{sing}). 
The point is here that Borcherds, while introducing his lift for arbitrary signature and also for general (scalar-valued) theta kernels, focuses on the Hermitian case. 
In particular, the geometric interpretation given in \cite{Bo2,Br1} only applies to this case.  
Hence from this aspect the lift $\Phi(z,f,\psi)$ and its features are new.
We put $\Lambda_{B,\psi}(z,f) = -d \Phi(z,f,\psi)$ and one sees that this again extends to a closed smooth $q$-form on $M$. Moreover, it essentially vanishes for $f \in M^!_k$.
One obtains as above (for a precise statement, see Theorem \ref{Main2}):
\begin{theorem}\label{IT4}
In the case of signature $(p,q)$, the lifts $\Lambda_{KM}$ and $\Lambda_{B,\psi}$ are adjoint.  
\end{theorem} 

This, almost axiomatic, introduction of the Borcherds lift for general signature does indeed lead to new information:  

\begin{theorem}\label{IT5}
Let $f \in H_k^+$. Then $\Lambda_{B,\psi}(z,f)$ is a harmonic representative of the Poincar\'e dual class of the cycle $Z(f)$. Moreover, we have the following equation of currents
\[
d[\Phi(z,f,\psi)] + \delta_{Z(f)} =[\Lambda_{B,\psi}(z,f)].
\]
In particular, the pair $\left(Z(f), \Phi(z,f,\psi)\right)$ defines a differential character in the sense of Cheeger and Simons (see \cite{C,CS}).
For $q=2$, we have in addition for the `classical' Borcherds lift
\[
dd^c[\Phi(z,f)] + \delta_{Z(f)} =[\Lambda_{B}(z,f)],
\]
i.e., $\Phi(z,f)$ is a Green's function for the divisor $Z(f)$.
\end{theorem} 

For $q=2$ the latter result already follows from \cite{Br1,Br2}.
We also briefly discuss the relationship of  $\Phi(z,f)$ to the Green's functions for the special divisors constructed by Oda and Tsuzuki \cite{OT} and Kudla \cite{K, Ku2}.




\smallskip

The results of this paper are subject to several extensions and generalizations. On one hand, one should be able to introduce suitable Borcherds lifts for unitary groups, since \cite{KM90} covers these groups as well. This is of particular interest, as the symmetric spaces are Hermitian in this case. In particular, this should lead to the explicit construction of Green's currents for cycles of complex codimension $q$. On the other hand,  Funke and Millson \cite{FM2} are currently developing a theory for special cycles with coefficients analogous to the Kudla-Millson lift. This also should give rise to new Borcherds lifts with geometric importance. 
On a more speculative note, it is a very interesting problem to generalize the present results  to special cycles of higher codimension. In the frame work of this paper, one would need to define a suitable Borcherds lift for the symplectic group $\Symp_n(\R)$. For this one would need an analogue of the space $H_k^+$ and henceforth of Theorem~\ref{IT1}. 
We hope to come back to these issues in the near future.

\smallskip

The paper is organized as follows. After setting up the basic notions in section 2, we discuss in section 3 the space of weak Maass forms $H_k$ in detail and prove Theorem~\ref{IT1}. In section 4, we consider the Fock model of the Weil representation and the Schwartz forms $\varphi_{KM}$ and $\psi$ in detail and derive Theorem~\ref{IT3}. We discuss the Kudla-Millson lift and the Borcherds lift in section 5, introducing for general signature the Borcherds lift $\Phi(z,f,\psi)$. The main result, Theorems~\ref{IT2} and \ref{IT4}, is proven in section 6. Finally, Theorem~\ref{IT5} is considered in section 7.

\smallskip

This project was first conceived and advanced during three visits from 1999 to 2002 of the second named author in Heidelberg. He would like to thank the \emph{Forschergruppe Arithmetik} at Heidelberg and E.~Freitag in particular for their generous hospitality during these visits.  He also thanks J.~Burgos and J.~Naranjo at the University of Barcelona for their patient hospitality in the summer 2002, where major work on this paper was done.
We also would like to thank J.~Burgos, E.~Freitag, W.~Kohnen, S.~Kudla, and J.~Millson for useful discussions and comments on this project.

%

\section{Basic Notations}
\label{sect2}

Let $V$ be a rational vector space over $\Q$ with a non-degenerate bilinear form $(\,,\,)$ of signature $(p,q)$; we assume $\dim V \geq 3$. Let $L$ be an even lattice in $V$ (that is $q(x) := \tfrac12(x,x) \in \Z$ for $x \in L$) and write $L^\#$ for its dual. We will write $V^-$ (respectively $L^-$) for the vector space $V$ ($\Z$-module $L$) together with the bilinear form $-(\,,\,)$.

We pick an orthogonal  basis $\{v_i\}$ of $V(\R)=V\otimes_\Q\R$ such that
$(v_{\alpha},v_{\alpha}) = 1$ for $\alpha = 1,\dots,p$ and
$(v_{\mu},v_{\mu}) = -1$ for $\mu = p+1,\dots,p+q$. The
corresponding coordinates we denote by $x_i$. Throughout the paper
we will use the subscript $\alpha$ for the `positive' variables
and $\mu$ for the `negatives' ones.
We realize the  symmetric space associated to $V$ as the set of negative $q$-planes in $V(\R)$:
\begin{equation*}
D \simeq \{ z \subset  V(\R);\; \dim z =q \quad \text{and} \quad
(\,,\,)|_z < 0\}.
\end{equation*} 
Occasionally, we will write
$D_{p,q}$ to emphasize the signature. The assignment $z \mapsto
z^{\perp}$ gives the identification  $D_{p,q} \simeq D_{q,p}$.
Putting $G = \SO_0(V(\R))$, the connected component of the
orthogonal group,  and letting $K$ be the maximal compact subgroup
of $G$ stabilizing $z_0 = \Span\{v_{\mu};\; p+1 \leq \mu \leq
p+q\}$, we certainly have $D \simeq G/K$.

For $z \in D$, we associate the standard majorant $(\,,\,)_z$ given by
\begin{equation*}
(x,x)_z = (x_{z^{\perp}},x_{z^{\perp}}) - (x_z,x_z),
\end{equation*}
where $x = x_z + x_{z^{\perp}} \in V(\R)$ is given by the
orthogonal decomposition $V(\R) = z \oplus z^{\perp}$.

Let $\mathfrak{g}$ be the Lie algebra of $G$ and $\mathfrak{g} = \mathfrak{p} + \mathfrak{k}$ its Cartan decomposition. Then $\mathfrak{p} \simeq \mathfrak{g}/  \mathfrak{k}$ is isomorphic to the tangent space at the base point of $D$, and with respect to the above basis of $V(\R)$ we have
\begin{equation}\label{mathfrakp}
\mathfrak{p} \simeq \left\{
\begin{pmatrix}
0 & X \\
{^tX} &0
\end{pmatrix}
;\,\, X \in M_{p,q}(\R) \right\} \simeq M_{p,q}(\R).
\end{equation}

For $q=2$, it is well known that $D$ is Hermitian, and we assume that the complex structure on $\mathfrak{p}$ is given by right multiplication with
$J = \left( \begin{smallmatrix}
 0&1 \\ -1&0 \end{smallmatrix} \right)$. 

Let $\G \subset G$ be a congruence subgroup of the orthogonal
group of $L$ fixing the discriminant lattice $L^{\#}/L$ 
and $X = \G \back D$ be the
associated locally symmetric space.
Throughout we assume that $\Gamma$ is torsion free such that $\Gamma$ acts freely on $D$, and $X$ is a real analytic manifold of dimension $pq$. Observe that for instance the principal congruence subgroup $\Gamma(N)$ of level $N$ (of the discriminant kernel of $\Orth(L)$), i.e., the the kernel of the natural homomorphism
$\Orth(L)\to \Orth\left((\tfrac{1}{N}L^\#)/L\right)$, is torsion free for $N\geq 3$ (which is seen similarly as Hilfssatz 6.5 in \cite{Fr} chapter II.6).  


For $x \in V(\R)$ with $q(x) >0$, we let
\begin{equation*}
D_x = \{ z \in D;\; z \perp x\}.
\end{equation*}
Note that $D_x$ is a subsymmetric space of type $D_{p-1,q}$ attached to the orthogonal group $G_x$, the pointwise stabilizer of $x$ in $G$. Put $\G_x = \G \cap G_x$. Then for $x \in L^{\#}$, the quotient
\begin{equation*}
Z(x) = \G_x \back D_x \longrightarrow X
\end{equation*}
defines a (in general relative) cycle in $X$. For $h \in
L^{\#}/L$ and $n\in \Q$, the group $\G$ acts on $L_{h,n} = \{x \in L +h;\; q(x)
= n \}$ with finitely many orbits, and we define the
composite cycle
\begin{equation*}
Z(h,n) = \sum_{ x \in \G \back L_{h,n}} Z(x).
\end{equation*}
Occasionally, we will identify  $Z(h,n)$ with its preimage in $D$.

Borcherds \cite{Bo2} and Bruinier \cite{Br1} (for $\Orth(2,p)$) use vectors of \emph{negative} length to define special cycles (which are divisors in that case); by switching to the space $V^-$ these are the same as the divisors for $D_{p,2}$ defined above by vectors of \emph{positive} length.

We orient $D$ and the cycles $D_x$ as in \cite{KM90}, page 130/131. Note that for $q$ even, $D_x$ and $D_{-x}$ have the same orientation while for $q$ odd, the opposite. Moreover, for $q=2$, this orientation coincides with the orientation given by the complex structure on $D$ and $D_x$.


\bigskip


Let $G'= \Mp_2(\R)$ be the two-fold cover of $\Sl_2(\R)$, realized
by the two choices of holomorphic square roots of $\tau \mapsto
j(g,\tau) = c\tau + d$; here $\tau \in \h=\{ w\in \C; \;
\Im(w)>0\}$, the upper half plane and $g = \left(
\begin{smallmatrix} a&b \\ c&d \end{smallmatrix} \right) \in
\Sl_2(\R)$. So elements of $\Mp_2(\R)$ are of the form $(g,
\phi(\tau))$ with $g \in \Sl_2(\R)$ and $\phi(\tau)$ a holomorphic
function such that $\phi(\tau)^2 = j(g,\tau)$. The multiplication
is given by $(g_1, \phi_1(\tau)) (g_2,\phi_2(\tau)) = (g_1g_2, \phi_1(g_2 \tau) \phi_2(\tau))$,
where $\Sl_2(\R)$ acts on $\h$ by linear fractional transformations.
For $z= r e^{i\theta} \in \C^\times$ with $\theta \in (-\pi,\pi]$
and $r$ positive, we take $\sqrt{z}= z^{1/2}= r^{1/2}e^{i\theta
/2}$. Occasionally we just write $g$ for $(g,\sqrt{j(g,\tau)})\in
G' $.

We write $K'$  for the inverse image of $\SO(2)\simeq \Uni(1)$  under the covering map $\Mp_2(\R) \longrightarrow \Sl_2(\R)$. Note that for $k_{\theta} \in \SO(2)$ with $k_{\theta} = \left(
\begin{smallmatrix}
\cos(\theta) & \sin(\theta) \\
-\sin(\theta) & \cos(\theta)
\end{smallmatrix} \right)$, ($\theta \in (-\pi,\pi]$), we obtain a character of $K'$ by the assignment
\begin{equation*}
\chi_{1/2}: (k_{\theta}, \pm \sqrt{j(k_{\theta},\tau)}) \mapsto \pm \sqrt{j(k_{\theta},i)}^{-1} = \pm e^{i  \theta/2}.
\end{equation*}

\smallskip

We denote by $\omega = \omega_V$ the Schr\"odinger model of the (restriction of the) Weil representation of $G'\times \Orth(V(\R))$ acting on $\mathcal{S}(V(\R))$, the space of Schwartz functions on $V(\R)$. We  have
\begin{equation*}
\omega(g) \varphi(x) = \varphi(g^{-1}x)
\end{equation*}
for $\varphi \in \mathcal{S}(V(\R))$ and $g \in \Orth(V(\R))$. The action of $G'$ is given as follows:
\begin{equation*}
\omega(m(a)) \varphi(x) = a^{(p+q)/2} \varphi(ax)
\end{equation*}
for $a>0$ and with $m(a) =  \left( \begin{smallmatrix} a&0\\0&a^{-1}
\end{smallmatrix} \right)$;
\begin{equation*}
\omega(n(b)) \varphi(x) = e^{\pi i b(x,x)} \varphi(x)
\end{equation*}
with $n(b) =  \left( \begin{smallmatrix} 1&b\\0&1
\end{smallmatrix} \right)$; and
\begin{equation*}
\omega(S) \varphi(x) = \sqrt{i}^{p-q} \hat{\varphi}(-x)
\end{equation*}
with $S =  \left( \left( \begin{smallmatrix} 0&-1\\1&0
\end{smallmatrix} \right), \sqrt{\tau} \right)$, and where $\hat{\varphi}(y) = \int_{V(\R)} \varphi(x) e^{2\pi i(x,y)} dx $ is the Fourier transform.
If $p+q$ is even, $(1,t) \in G'$ acts trivially, otherwise by
multiplication by $t$.

For $\varphi \in \mathcal{S}(V(\R))$ and $h \in L^{\#}/L$, we define the theta function
\begin{equation*}
\theta(g',\varphi,h) = \sum_{\la \in L +h} \omega(g') \varphi(\la) \qquad \qquad (g' \in G').
\end{equation*}

We write $\G'$ for the inverse image of $\Sl_2(\Z)$ in $\Mp_2(\R)$.
For the generators $S$ and $T=(n(1),1)$ of $\G'$, we then have
\begin{equation}\label{weilt}
\theta(Tg',\varphi,h) = e^{\pi i (h,h)} \theta(g',\varphi,h),
\end{equation}
and by Poisson summation
\begin{equation}\label{weils}
\theta(Sg',\varphi,h) = \frac{\sqrt{i}^{p-q}}{\sqrt{|L^{\#}/L|}} \sum_{h'\in L^{\#}/L} e^{- 2 \pi i (h,h')} \theta(g',\varphi,h').
\end{equation}
The equations \eqref{weilt} and \eqref{weils} define a 
representation $\varrho_L$ of $\G'$ acting on the group algebra
$\C[L^{\#}/L]$, whose standard basis elements we denote by
$\frake_{h}$, $h \in {L^{\#}/L}$. Defining the vector
\begin{equation*}
\Theta(g',\varphi,L) :=  (\theta(g',\varphi,h))_{h \in {L^{\#}/L}} = \sum_{h \in {L^{\#}/L}} \theta(g',\varphi,h)\frake_h,
\end{equation*}
we then have
\begin{equation*}
\Theta(\g g',\varphi,L) = \varrho_L(\g) \Theta(g',\varphi,L)
\end{equation*}
for all $\g \in \G'$. 
Note that with respect to the standard scalar product $\langle\,,\,\rangle $ on  $\C[L^{\#}/L]$ (linear in the first and anti-linear in the second variable),  we have that $\varrho_{L^-} = \varrho_L^{\ast} = \bar{\varrho_L}$;
i.e., $L$ and $L^-$ give rise to dual representations. From
\eqref{weilt} and \eqref{weils} it is clear that $\varrho_L$
coincides with the representation   $\varrho_L$ considered in
\cite{Bo1} and \cite{Br1}. Ultimately, $\varrho_L$ goes back to
Shintani \cite{Sh}, which is also a good reference for the above
discussion.

\smallskip

Let $\varphi$ now be an eigenfunction under the action of $K'$; i.e., $\omega(k_{\theta})\varphi = \chi^r_{1/2}(k_{\theta}) \varphi$ for some $r \in \Z$, so that
\begin{equation*}
\Theta(g'k_{\theta},\varphi,L) = \chi^r_{1/2}(k_{\theta})\Theta(g',\varphi,L).
\end{equation*}
Then we can associate to $\Theta(g',\varphi,L)$ a (vector-valued) function on the upper half plane in the usual way: We let $g'_{\tau} =
\left( \begin{smallmatrix}
1&u \\0& 1
\end{smallmatrix} \right)
\left( \begin{smallmatrix}
v^{1/2}&0 \\0& v^{-1/2}
\end{smallmatrix} \right)$
with $\tau = u+iv \in \h$ be the standard element moving the base point $i \in \h$ to $\tau$ and define
\begin{equation*}
\Theta(\tau,\varphi,L) := j(g'_{\tau},i)^{r/2} \Theta(g'_{\tau},\varphi,L).
                       = \sum_{h \in L^{\#}/L} \sum_{\la \in L +h} \varphi(\la,\tau,z) \frake_h
\end{equation*}
with 
\[
\varphi(\la,\tau,z) =  j(g'_{\tau},i)^{r/2} \omega(g'_{\tau})\varphi(\la).
\]
Hence, for $(\g,\phi) \in \G'$,
\begin{equation}\label{trafo}
\Theta(\g\tau,\varphi,L) = \phi(\tau)^r \varrho_L(\g,\phi) \Theta(\tau,\varphi,L);
\end{equation}
i.e., $\Theta(\tau,\varphi,L)$ is a $C^{\infty}$-automorphic form of weight $r/2$ with respect to the representation $\varrho_L$. (Note however, that it is usually not an eigenfunction of the Laplacian).

We denote the real analytic functions on $\h$ satisfying the transformation property (\ref{trafo}) with weight $k$ by $A_{k,L}$.



\smallskip


From now on we will frequently drop the lattice $L$ from the argument of $\Theta(g',\varphi,L)$.

\smallskip


The space of $K$-invariant Schwartz functions $\mathcal{S}(V(\R))^K$ is of particular interest as we have
\begin{equation*}
\mathcal{S}(V(\R))^K \simeq [\mathcal{S}(V(\R)) \otimes C^{\infty}(D)]^G,
\end{equation*}
where the isomorphism is giving by evaluation at the base point
$z_0$ of $D$.  Note that for such $\varphi(x,z) \in [\mathcal{S}(V(\R)) \otimes C^{\infty}(D)]^G$, the theta function $\Theta(g',z,\varphi)$ is $\G$-invariant as a function of $z$; hence, it descends to a function on $X$.

The Gaussian on $V(\R)$ is given by $\varphi_0(x) = e^{- \pi \sum_{i=1}^{p+q} x_i^2}$. Certainly $\varphi_0 \in \mathcal{S}(V(\R))^{K}$, and the corresponding function in
$[\mathcal{S}(V(\R)) \otimes C^{\infty}(D)]^{ G}$ is given by
\begin{equation*}
\varphi_0(x,z) = e^{-\pi (x,x)_z}.
\end{equation*}
Occasionally we will write $\varphi_0^{p,q}$ to emphasize the signature. We have
\begin{equation*}
\Theta(\tau,z,\varphi_0) = \sum_{h \in L^{\#}/L} 
\sum_{\la \in L +h} \varphi_0(\la,\tau,z) \frake_h \in A_{(p-q)/2,L} \otimes
C^{\infty}(X),
\end{equation*}
and
\begin{equation*}
\varphi_0(\la,\tau,z) = v^{q/2} \exp\big( \pi i((\la,\la)u + (\la,\la)_ziv)\big) = \exp\big(\pi i( (\la_{z^{\perp}},\la_{z^{\perp}})\tau + (\la_z,\la_z)\bar{\tau})\big).
\end{equation*}

\section{Weak Maass forms}
\label{sectmaass}

In this section we introduce a new space of Maass wave forms. In particular, we establish a pairing with holomorphic modular forms and obtain a duality theorem for this pairing.
The results of this section can be either viewed as an analytic version of the Serre duality result in \cite{Bo2} section 3, or as an algebraic approach to 
\cite{Br1} chapter 1. 

Let $k\in \frac{1}{2}\Z$ with $k\neq 1$. 
(In our later applications $k$ will be smaller than $1$. 
However, we do not need that here.) 
Moreover, let $\Gamma''\leq \Gamma'$ be a subgroup of finite index.

We write $H_{k,L}(\Gamma'')$ for the space  of {\em weak Maass forms} of weight $k$ with representation $\varrho_L$ for the group $\Gamma''$. By definition, this is the space of twice continuously differentiable functions $f:\H\to\C[L^\#/L]$ satisfying:

\begin{enumerate}
\renewcommand{\labelenumi}{\roman{enumi})}    
\item
$f(\gamma\tau) = \phi(\tau)^{2k} \varrho_L(\gamma,\phi) f(\tau)$ for all $(\gamma,\phi)\in \Gamma''$;
\item 
there is a $C>0$ such that for any cusp $s\in \Q \cup \{\infty\}$ of $\Gamma''$ and $(\delta,\phi)\in \Gamma'$ with $\delta\infty=s$ the function 
$f_s(\tau) = \phi(\tau)^{-2k} \varrho_L^{-1}(\delta,\phi) f(\delta\tau)$ satisfies
$f_s(\tau)=O(e^{C v})$ as $v\to \infty$ (uniformly in $u$, where $\tau=u+iv$);
\item 
$\Delta_k f(\tau)=0$, 
where 
\begin{equation}\label{defdelta}
\Delta_k = -v^2\left( \frac{\partial^2}{\partial u^2}+ \frac{\partial^2}{\partial v^2}\right) + ikv\left( \frac{\partial}{\partial u}+i \frac{\partial}{\partial v}\right)
\end{equation}
denotes the usual hyperbolic Laplace operator in weight $k$.
\end{enumerate} 


Since $\Delta_k$ is an elliptic differential operator, such a function $f$ will automatically be real analytic.
We will mainly work with the full group $\Gamma'$ and therefore abbreviate $H_{k,L}=H_{k,L}(\Gamma')$.
The transformation property (i) implies that any $f\in H_{k,L}$ has a Fourier expansion 
\[
f(\tau)= \sum_{h\in L^\#/L}\sum_{n\in \Q} a(h,n;v)e(nu)\frake_h,
\]
where $e(u)=e^{2\pi i u}$, as usual. The coefficients $a(h,n;v)$ vanish unless $n-q(h)\in \Z$. In particular, the denominators of the indices $n$ of all non-zero coefficients $a(h,n;v)$ are bounded by the level $N$ of the lattice $L$.
Because of property (iii), the coefficients $a(h,n;v)$ satisfy the second order differential equation $\Delta_k a(h,n;v)e(nu)=0$ as functions in $v$. 
If $n=0$, one finds that $a(h,0;v)$ is a linear combination of  $1$ and  $v^{1-k}$.
If $n\neq 0$, then, writing $a(h,n;v)$ as $b(2\pi n v)$, it is easily seen that $b(w)$ is a solution of the second order linear differential equation
\begin{align*}
\frac{\partial^2}{\partial w^2} b(w)- b(w) + \frac{k}{w}\left( 
\frac{\partial}{\partial w} b(w) +b(w)\right)=0,
\end{align*} 
which is independent of $h$ and $n$.
It is immediately checked that $e^{-w}$ is a solution. A second, linearly independent solution is found by reduction of the order. Here we choose the function
\[
H(w)=e^{-w}\int_{-2w}^\infty e^{-t}t^{-k}\, dt.
\] 
The integral converges for $k<1$ and can be holomorphically continued in $k$ (for $w\neq 0$) in the same way as the Gamma function.
If $w<0$, then $H(w)=e^{-w}\Gamma(1-k,-2w)$, where $\Gamma(a,x)$ denotes the incomplete Gamma function as in \cite{AS} (6.5.3). 
The function $H(w)$ has the asymptotic behavior
\[
H(w)\sim \begin{cases}  
(2|w|)^{-k} e^{-|w|},&\text{for $w\to -\infty$,}\\
(-2w)^{-k} e^{w},&\text{for $w\to +\infty$.}
\end{cases}
\]
(Clearly the functions $|w|^{k/2}e^{-w/2}$ and $|w|^{k/2}H(w/2)$ are the special values of the standard Whittaker functions $W_{\nu,\mu}(|w|)$ and $M_{\nu,\mu}(|w|)$ for $\nu =\sgn(n)k/2$, $\mu=k/2-1/2$ as in \cite{AS} chapter 13.)
We find that
\[
a(h,n;v)=\begin{cases}
a^+(h,0)+a^-(h,0)v^{1-k},& \text{if $n=0$,}\\
a^+(h,n) e^{-2\pi nv}+a^-(h,n) H(2\pi nv),&\text{if $n\neq 0$,}
\end{cases}
\]
with complex coefficients $a^\pm (h,n)$. Thus any weak Maass form $f$ of  weight $k$ has a unique decomposition $f=f^++f^-$, where
\begin{subequations}
\label{deff}
\begin{align}
\label{deff+}
f^+(\tau)&= \sum_{h\in L^\#/L}\sum_{n\in \Q} a^+(h,n) e(n\tau)\frake_h,\\
\label{deff-}
f^-(\tau)&= \sum_{h\in L^\#/L} \bigg(a^-(h,0) v^{1-k}
+\sum_{\substack{n\in \Q\\ n\neq 0}} a^-(h,n) H(2\pi nv) e(nu)\bigg)\frake_h.
\end{align}
\end{subequations}
Note that if $f$ satisfies condition (ii) above, then all but finitely many $a^+(h,n)$ (respectively  $a^-(h,n)$) with negative (respectively positive) index $n$ vanish.  

Let us briefly recall the Maass raising and lowering operators on non-holomorphic modular forms of weight $k$.
They are defined as the differential operators
\begin{align*}
R_k  =2i\frac{\partial}{\partial\tau} + k v^{-1} \qquad \text{and} \qquad 
L_k  = -2i v^2 \frac{\partial}{\partial\bar{\tau}}.
\end{align*}
The raising operator $R_k$ maps $A_{k,L}$ to $A_{k+2,L}$, and the lowering operator $L_k$  maps $A_{k,L}$ to $A_{k-2,L}$.
The Laplacian $\Delta_k$ can be expressed in terms of $R_k$ and $L_k$ by
\begin{equation}\label{deltalr}
-\Delta_k = L_{k+2} R_k +k = R_{k-2} L_{k}.
\end{equation}

The following lemma is proved by a straightforward computation.

\begin{lemma}\label{fourier}
Let $f\in H_{k,L}$ be a weak Maass form  of weight $k$ and write $f=f^+ + f^-$ as in \eqref{deff}. Then  
\[
L_k f = L_k f^- = -2 v^{2-k}\sum_{h\in L^\#/L}\bigg(  (k-1) a^-(h,0)+ \sum_{\substack{n\in \Q\\ n\neq 0}} a^-(h,n) (-4\pi n)^{1-k}e(n\bar\tau)\bigg)\frake_h.
\]
\end{lemma}

We write $M^!_{k,L}(\Gamma'')$ for the space of holomorphic $\C[L^\#/L]$-valued functions on $\H$ satisfying the transformation property (i) above and being meromorphic at the cusps of $\Gamma''$. We call such modular forms {\em weakly holomorphic}.
Identity \eqref{deltalr} implies that $M^!_{k,L}(\Gamma'')\subset H_{k,L}(\Gamma'')$.
If we also require holomorphicity (vanishing) at the cusp, we obtain the space of holomorphic modular forms $M_{k,L}(\Gamma'')$ (cusp forms $S_{k,L}(\Gamma'')$). Finally, if $\Gamma''$ is the full modular group $\Gamma'$, then we will briefly write $M^!_{k,L}$, $M_{k,L}$, and $S_{k,L}$ for these spaces of modular forms.   
They were first considered by Borcherds \cite{Bo1} and later in  \cite{Br1}. 

The lattice $L^-$ gives rise to the dual representation $\varrho_L^{\ast}= \bar{\varrho_L}$. Hence, we can consider $A_{k,L^-}$, etc. also as the space of modular forms with respect to the dual representation $\varrho_L^{\ast}$ of $\varrho_L$. We will make frequent use of this fact.

Recall that the Petersson scalar product on $M_{k,L}$ is  given by
\begin{equation}\label{pet}
(f,g)_{k,L} = \int_{\Gamma' \back \H} \langle f,g \rangle v^{k}\,
d\mu
\end{equation}
for $f,g\in M_{k,L}$, whenever the integral converges absolutely. 
Here $d\mu=\frac{du\,dv}{v^2}$ denotes the usual invariant volume form on $\H$.

\begin{proposition}\label{defmap}
The assignment 
$f(\tau)\mapsto \xi_k(f)(\tau):=v^{k-2} \overline{L_k f(\tau)} = R_{-k} v^k\overline{ f(\tau)}$ defines an antilinear mapping 
\begin{equation}\label{defxi}
\xi_k:H_{k,L}\longrightarrow M^!_{2-k,L^-}.
\end{equation}
Its kernel is  $M^!_{k,L}\subset H_{k,L}$.
\end{proposition}

\begin{proof}
By Lemma \ref{fourier}, the function $\xi_k(f)(\tau)$ is holomorphic on $\H$.  It vanishes if and only if $f^-$ vanishes. It is meromorphic at the cusp via the growth condition on $f$ . The transformation behavior of $\xi_k(f)(\tau)$ is easily checked. 
\end{proof}


We denote the inverse image of the space of holomorphic cusp forms $S_{2-k,L^-}$ under the mapping $\xi_k$  by $H^+_{k,L}$. Hence, if $f\in H^+_{k,L}$, then the Fourier coefficients $a^-(h,n)$ with non-negative index $n$ vanish, so $f^-$ is rapidly decreasing for $v \to \infty$. Clearly $M^!_{k,L}\subset H^+_{k,L}$. 
Moreover, if $k\geq 2$, then  $f^-$ vanishes for all $f\in H^+_{k,L}$.

Let now $f\in H_{k,L}$ and write its Fourier expansion as in \eqref{deff}. Then we call the Fourier polynomial 
\begin{align}\label{principal}
P(f)(\tau)&= \sum_{h\in L^\#/L}\sum_{\substack{n\in \Q\\n\leq 0}} a^+(h,n) e(n\tau)\frake_h
\end{align} 
the {\em principal part} of $f$. Observe that if $f\in H^+_{k,L}$, then $f-P(f)$ is exponentially decreasing as $v\to\infty$. (This property could actually be used to define the space $H^+_{k,L}$ in an alternative way.)

\begin{lemma}\label{growth1}
If $f\in H_{k,L}$, then there is a constant $C>0$ such that $f(\tau)=O(e^{C/v})$ as $v\to 0$, uniformly in $u$.
\end{lemma}

\begin{proof}
This follows from the transformation behavior and the growth of $f$ in the same way as the analogous statement for holomorphic modular forms.
\end{proof}

Later, we will need the following growth estimate for the Fourier coefficients of weak Maass forms.

\begin{lemma}\label{growth2}
Let $f\in H_{k,L}$ and write its Fourier expansion as in \eqref{deff}.
Then there is a constant $C>0$ such that the Fourier coefficients satisfy
\begin{align*}
a^+(h,n)&= O\left(e^{C\sqrt{|n|}}\right),\quad n\to +\infty,\\
a^-(h,n)&= O\left(e^{C\sqrt{|n|}}\right),\quad n\to -\infty.
\end{align*}
If $f\in H^+_{k,L}$, then the $a^-(h,n)$ actually satisfy the stronger bound
$a^-(h,n)=O(|n|^{k/2})$ as $n\to-\infty$.
\end{lemma}

\begin{proof}
To prove the asymptotic for the $a^-(h,n)$ we consider the weakly holomorphic modular form $\xi_k(f)\in M^!_{2-k,L^-}$.
By Lemma \ref{fourier} and the formula for the Fourier coefficients we have
\begin{equation}\label{growth21}
2 a^-(h,n)(-4\pi n)^{1-k} = -\int_{0}^1  \left\langle v^{k-2}\overline{L_k f(\tau)},\frake_h \right\rangle e(n\tau)\, du.
\end{equation}
Thus, according to Lemma \ref{growth1}, we get
\[
a^-(h,n)\ll |n|^{k-1} \int_{0}^1 e^{C/v} e^{-2\pi nv}\, du
\]
for all positive $0<v\leq 1$ with some positive constant $C$ (independent of $v$ and $n$).
If we take $v$ equal to $1/\sqrt{|n|}$, we see that
\[
a^-(h,n)\ll |n|^{k-1} e^{C\sqrt{|n|}} e^{2\pi \sqrt{|n|}}
\]
for all $n<0$, proving the first assertion on the $a^-(h,n)$. 

From $\xi_k(f)\in M^!_{2-k,L^-}$ it can be deduced that the individual functions $f^+$ and $f^-$ in the splitting $f=f^++f^-$ 
also satisfy the estimate of Lemma \ref{growth1}. 
We may apply the above argument to 
\begin{equation}\label{growth22}
 a^+(h,n) = \int_{0}^1  \left\langle f^+(\tau),\frake_h \right\rangle e(-n\tau)\, du
\end{equation}
to derive the estimate for the $a^+(n,h)$ as $n\to +\infty$.

If $f\in H^+_{k,L}$, then $\xi_k(f)\in S_{2-k,L^-}$ is a holomorphic cusp form. Hence the usual Hecke bound for the Fourier coefficients of cusp forms implies that the left hand side of \eqref{growth21} is bounded by some constant times $|n|^{1-k/2}$  for all $n<0$.
Thus $a^-(h,n)=O(|n|^{k/2})$ as $n\to -\infty$.
\end{proof}

The estimates of Lemma \ref{growth2} are far from being optimal. However, they are sufficient for our purposes. Observe that all the above results have obvious generalizations to finite index subgroups $\Gamma''\leq \Gamma'$.

\bigskip

Put $\kappa=2-k$. We now define a bilinear pairing between the spaces $M_{\kappa,L^-}$  and $H^+_{k,L}$ by putting 
\begin{equation}\label{defpair}
\{g,f\}=\big( g,\, \xi_k(f)\big)_{\kappa,L^-} 
\end{equation}
for $g\in M_{\kappa,L^-}$ and $f\in H^+_{k,L}$. 

\begin{proposition}\label{proppair}
Let $g\in M_{\kappa,L^-}$ with Fourier expansion $g=\sum_{h,n}b(h,n) e(n\tau)\frake_h$, and $f\in H^+_{k,L}$ with Fourier expansion as in \eqref{deff}.
Then the pairing \eqref{defpair} of $g$ and $f$ is determined by the principal part of $f$. It is equal to
\begin{equation}\label{pairalt}
\{g,f\}= \sum_{h\in L^\#/L} \sum_{n\leq 0}  a^+(h,n) b(h,-n).
\end{equation}
\end{proposition}

\begin{proof}
We begin by noticing that $\langle g,\overline{f}\rangle d\tau$ is a $\Gamma'$-invariant $1$-form on $\H$.
We have
\begin{align*}
d\left( \langle g,\overline{f}\rangle d\tau\right) = \bar\partial \left(\langle g,\overline{f}\rangle d\tau\right)
= \left\langle g,\overline{\tfrac{\partial}{\partial\bar \tau} f} \right\rangle d\bar\tau\, d\tau
= -\langle g,\overline{L_k f} \rangle \,d\mu.
\end{align*}
Hence, by Stokes' theorem we get
\begin{align*}
 \int_{\calF_t} \langle g,\overline{L_k f} \rangle \,
d\mu
=-\int_{\partial\calF_t}\langle g,\overline{f}\rangle\, d\tau
=\int_{-1/2}^{1/2} \langle g(u+it),\overline{f(u+it)}\rangle\, d u,
\end{align*}
where 
\begin{equation}\label{Fdomain}
\calF_t=\{ \tau\in \H;\; \text{$|\tau|\geq 1$, $|u|\leq 1/2$ , and $v\leq t$}\}
\end{equation}
denotes the truncated fundamental domain for the action of $\Sl_2(\Z)$ on $\H$.
In the last line we have used the invariance of $\langle g,\overline{f}\rangle d\tau$. If we insert the Fourier expansions of $f$ and $g$, we see that the integral over $u$ picks out the $0$-th Fourier coefficient of $\langle g,\overline{f}\rangle$. Thus  
\begin{align*}
\int_{\calF_t} \langle g,\overline{L_k f} \rangle \,
d\mu
&=\sum_h\sum_{n\leq 0} a^+(h,n)b(h,-n)+  O(e^{-\eps t}) 
\end{align*}
for some $\eps>0$. We finally obtain
\begin{align*}
\{g,f\}&= \lim_{t\to\infty} \int_{\calF_t} \langle g,\overline{L_k f} \rangle \,
d\mu =\sum_{h\in L^\#/L} \sum_{n\leq 0}  a^+(h,n) b(h,-n)
\end{align*}
as asserted.
\end{proof}

Observe that the definition of the pairing \eqref{defpair} immediately implies that $\{g,f\}=0$ for all  $f\in M^!_{k,L}$. By \eqref{pairalt} we get nontrivial relations among the coefficients of modular forms in $M_{\kappa,L^-}$. These are also easily obtained by means of the residue theorem on the Riemann sphere.
Moreover, it is clear that $\{g,f\}=0$ for all Eisenstein series $g\in M_{\kappa,L^-}$.  

\begin{theorem}\label{thmpair}
The pairing between the quotient $H^+_{k,L}/M^!_{k,L}$ and $S_{\kappa,L^-}$ induced by \eqref{defpair} is non-degenerate.  
\end{theorem}

It suffices to show that the mapping $\xi_k:H^+_{k,L}\to S_{\kappa,L^-}$ is surjective.
This is an immediate consequence of

\begin{theorem}\label{thmsurj}
The mapping $\xi_k:H_{k,L}\to M^!_{\kappa,L^-}$ defined in Proposition \ref{defmap} is surjective. 
\end{theorem}

Before beginning with the proof, we need to introduce some notation.

Let $\Gamma''\leq \Gamma'$ be a normal subgroup of finite index. 
We write $X(\Gamma'')$ for the compact modular curve corresponding to $\Gamma''$  and $\pi:\H\to X(\Gamma'')$ for the canonical map.
%
%
The Poincar\'e metric $1/v^2$ on $\H$ induces a K\"ahler metric on $X(\Gamma'')$ (with logarithmic singularities at the cusps). We write $*$ for the corresponding Hodge star operator.

We denote the sheaf of holomorphic functions (respectively $1$-forms) on $X(\Gamma'')$ by $\calO$ (respectively $\Omega$) and the sheaf of $C^\infty$ differential forms of type $(p,q)$ by $\calE^{p,q}$.
If $D$ is a divisor on $X(\Gamma'')$, then we write $\calO_D$ for the sheaf corresponding to $D$. The sections of $\calO_D$ over an open subset $U\subset X(\Gamma'')$ are given by meromorphic functions $f$ satisfying $\dv(f)\geq -D$ on $U$.

If $\calL$ is any $\calO$-module on $X(\Gamma'')$, we denote the $\calL$-valued $C^\infty$ differential forms on an open set $U\subset X(\Gamma'')$, that is, the sections of $\calE^{p,q}\otimes_\calO \calL$ over $U$, by $\calE^{p,q}(U,\calL)$.
Moreover, we write $\calL_{k,L}$ for the $\calO$-module sheaf of modular forms of weight $k$ with representation $\varrho_L$ on $X(\Gamma'')$. If $U\subset X(\Gamma'')$ is open, then the sections $\calL_{k,L}(U)$ are holomorphic $\C[L^\#/L]$-valued functions on the open subset $\pi^{-1}(U)\subset\H$ satisfying the transformation law of modular forms of weight $k$ with representation $\varrho_L$ for $\Gamma''$, and being holomorphic at the cusps.  If $\Gamma''$ acts freely on $\H$, then $\calL_{k,L}$ is a holomorphic vector bundle.
A Hermitean metric (with logarithmic singularities at the cusps) on it is given by the Hermitean scalar product $(f,g)_\tau=\langle f,g\rangle v^k$ on the fiber over $\pi(\tau)$, where $\tau=u+iv \in \H$ and $\langle f,g\rangle$ denotes the standard scalar product on $\C[L^\#/L]$ as before.

The dual vector bundle $\calL_{k,L}^*$ of $\calL_{k,L}$ can be identified with the vector bundle $\calL_{-k,L^-}$ of modular forms of weight $-k$ with representation $\varrho_L^*$ on $X(\Gamma'')$. The mapping $f\mapsto v^k \bar f$ defines an anti-linear bundle isomorphism $\calL_{k,L}\to \calL_{k,L}^*$. It induces a Hodge star operator 
\[
\bar{*}_{\calL}: \calE^{p,q}\otimes_\calO \calL_{k,L}\longrightarrow \calE^{1-p,1-q}\otimes_\calO \calL_{k,L}^*,\qquad 
\bar{*}_{\calL}(\phi\otimes f)=(*\bar \phi) \otimes (v^k \bar f), 
\]
on $\calL_{k,L}$-valued $C^\infty$ differential forms on $X(\Gamma'')$ (see \cite{We} chapter V.2). 
It is easily verified that
\begin{subequations}
\label{star}
\begin{align}
\bar{*}_{\calL}(f)&=v^{k-2}\bar f \,du\,dv, \qquad\text{for $f \in \calE^{0,0}(U,\calL_{k,L})$,}\\
\bar{*}_{\calL}(f d\tau)&=i v^k \bar f d\bar\tau,\qquad 
\text{for $f d\tau \in \calE^{1,0}(U,\calL_{k,L})$,}\\
\bar{*}_{\calL}(\bar f d\bar \tau)&=i v^k  f d \tau,\qquad 
\text{for $\bar f d\bar\tau \in \calE^{0,1}(U,\calL_{k,L})$.}
\end{align}
\end{subequations}
The Laplace operator on differential forms in $\calE^{p,q}(U,\calL_{k,L})$ is given by
\[
\overline{\Box}=\bar \partial\bar{*}_{\calL} \bar\partial\bar{*}_{\calL}+\bar{*}_{\calL}\bar \partial\bar{*}_{\calL} \bar\partial.
\]
By a straightforward computation using \eqref{star} it can be shown that 
\begin{align}\label{laplcompl}
\overline{\Box}f=\bar{*}_{\calL}\bar \partial\bar{*}_{\calL} \bar\partial \,f
=R_{k-2}L_k f =-\Delta_k f
\end{align}
for functions $f\in \calE^{0,0}(U,\calL_{k,L})$.

\begin{proof}[Proof of Theorem \ref{thmsurj}]
Let $\Gamma''\leq \Gamma'$ be a normal subgroup of finite index that acts freely on $\H$.
We first prove the analogous result for the group $\Gamma''$, namely that the mapping $\xi_k:H_{k,L}(\Gamma'')\to M^!_{\kappa,L^-}(\Gamma'')$ is surjective.

We consider the Riemann surface $X=X(\Gamma'')$ using the above notation. Let $s_1, \dots, s_r$ be the cusps of $X$ and write $D=\sum_i s_i$ for the divisor on $X$ given by the cusps.
Let $n$ be a positive integer.
By tensoring the Dolbeault resolution of the structure sheaf $\calO$ with the locally free $\calO$-module $\calL_{k,L}\otimes_\calO \calO_{nD}$, we get the exact sheaf sequence
\[
\xymatrix{
0\ar[r]& \calO \otimes \calL_{k,L}\otimes \calO_{nD}\ar[r] &\calE^{0,0} \otimes \calL_{k,L}\otimes \calO_{nD} \ar[r]^{\bar\partial\otimes 1\otimes 1} &\calE^{0,1} \otimes \calL_{k,L}\otimes \calO_{nD} \ar[r]& 0.
}
\]
Since $\calE^{p,q} \otimes \calL_{k,L}\otimes \calO_{nD}$ is a fine sheaf, it is acyclic. Hence we obtain the following long exact cohomology sequence:
\begin{align}\label{cohsequ}
0\longrightarrow  (\calL_{k,L}\otimes \calO_{nD})(X)& \longrightarrow \calE^{0,0}(X, \calL_{k,L}\otimes \calO_{nD}) \\
\nonumber
&\longrightarrow \calE^{0,1}(X,\calL_{k,L}\otimes \calO_{nD}) 
\longrightarrow H^1(X,\calL_{k,L}\otimes \calO_{nD}) \longrightarrow 0.
\end{align}

We claim that $H^1(X,\calL_{k,L}\otimes \calO_{nD})$ vanishes, if $n$ is large.
In fact, using Serre duality (see \cite{We} chapter V Theorem 2.7) we find
\[
H^1(X,\calL_{k,L}\otimes \calO_{nD})\cong H^0(X,\Omega\otimes\calL_{k,L}^*\otimes \calO_{nD}^*)
\cong H^0(X,\calL_{\kappa,L^-}\otimes \calO_{-nD}).
\]
Because the number of zeros (counted with multiplicities) of holomorphic modular forms of fixed weight $\kappa$ is bounded, the latter cohomology group vanishes if $n$ is sufficiently large. This proves the claim.


Now let $g\in M^!_{\kappa,L^-}(\Gamma'')$. We want to show that there is an $f\in H_{k,L}(\Gamma'')$ such that $\xi_k(f)=g$.
We chose a positive integer $n$ greater or equal to the orders of the poles of $g$ at the cusps and such that $H^1(X,\calL_{k,L}\otimes \calO_{nD})$ vanishes.
Then the $1$-form $g \,d \tau$ defines a global holomorphic section in $\calE^{1,0}(X,\calL_{k,L}^*\otimes \calO_{nD})$. Applying the Hodge star operator, we get the $1$-form
$iv^{-k}\bar g\, d\bar\tau\in \calE^{0,1}(X,\calL_{k,L}\otimes \calO_{nD})$.  
By virtue of the exact sequence \eqref{cohsequ} we find that there is a function $f\in \calE^{0,0}(X, \calL_{k,L}\otimes \calO_{nD})$ satisfying
$\bar\partial f=iv^{-k}\bar g\, d\bar\tau$.
But this is equivalent to saying
\[
-v^{k-2} L_k f\, d\bar\tau= \bar g\, d\bar\tau.
\]
We are left with showing that $\Delta_k f=0$. This follows from
\[
- \Delta_k f =  \overline{\Box} f
= \bar{*}_{\calL}\bar \partial\bar{*}_{\calL} \bar\partial \,f
= \bar{*}_{\calL}\bar \partial\bar{*}_{\calL} \, iv^{-k}\bar g\, d\bar\tau
= -\bar{*}_{\calL}\bar \partial \, g d \tau
=0
\]
completing the proof of Theorem \ref{thmsurj} for the group $\Gamma''$.

For the full modular group $\Gamma'$ the assertion follows by considering the $\Gamma'/\Gamma''$-invariant subspaces in \eqref{cohsequ}, using the fact that the cohomology of a finite group acting on a $\C$-vector space vanishes.
\end{proof}

\begin{corollary}\label{corsurj}
The following sequences are exact:
\begin{gather*}
\xymatrix{
0\ar[r]& M_{k,L}^! \ar[r]& H_{k,L} \ar[r]^{\xi_k}&  M^!_{\kappa,L^-} \ar[r] & 0,
}\\
\xymatrix{
0\ar[r]& M_{k,L}^! \ar[r]& H_{k,L}^+ \ar[r]^{\xi_k}&  S_{\kappa,L^-} \ar[r] & 0.
}
\end{gather*}
\end{corollary}


We will also need a slight variation of the above duality result.  To formulate it, we introduce a second pairing between the spaces $H^+_{k,L}$ and $M_{\kappa,L^-}$. If  $g\in M_{\kappa,L^-}$ with Fourier expansion $g=\sum_{h,n}b(h,n) e(n\tau)\frake_h $, and $f\in H^+_{k,L}$ with Fourier expansion as in \eqref{deff}, then we put
\begin{equation}\label{defpair'}
\{g,f\}'= \sum_{h\in L^\#/L} \sum_{n< 0}  a^+(h,n) b(h,-n).
\end{equation}
This definition differs from \eqref{pairalt} by  the fact that we only sum over negative $n$.
We let $M^{!!}_{k,L}$ be the subspace of $M^!_{k,L}$ consisting of those weakly holomorphic modular forms $f$ whose constant term $\sum_h a^+(h,0)\frake_h $ in the Fourier expansion is orthogonal to the constant terms of all $g\in M_{\kappa,L^-}$ with respect to $\langle\cdot,\cdot\rangle$.
If $k<0$ or if $k=0$ and $\varrho_L$ does not contain the trivial representation, then, using Eisenstein series in $M_{\kappa,L^-}$, one sees that $M^{!!}_{k,L}$ simply consists of those $f\in M^!_{k,L}$ with vanishing constant term. However, otherwise the fact that Eisenstein series of weight $\leq 2$ may be non-holomorphic implies that some $f\in M^{!!}_{k,L}$ may have non-vanishing constant term. 
%
%
If $f\in M^{!!}_{k,L}$ or $g\in S_{\kappa,L^-}$, then $\{g,f\}=\{g,f\}'$. In particular the pairing $\{g,f\}'$ vanishes, if $f\in M^{!!}_{k,L}$.

\begin{corollary}\label{corpair}
The pairings between $H^+_{k,L}/M^{!!}_{k,L}$ and $M_{\kappa,L^-}$ (respectively $H^+_{k,L}/M^{!}_{k,L}$ and $S_{\kappa,L^-}$) induced by \eqref{defpair'} are non-degenerate.  
\end{corollary}

This can be proved using Theorem \ref{thmpair} and Eisenstein series in $M_{\kappa,L^-}$. The details are  left to the reader.

\begin{remark}\label{Hejhal}
If $h\in L^\#/L$ and $m\in \Z+q(h)$ is negative, the Hejhal Poincar\'e series $F_{h,m}(\tau,\kappa/2)$ defined in \cite{Br1} Definition 1.8 and Proposition 1.10 (see also \cite{He} and \cite{Ni}) are examples of weak 
Maass forms in $H^+_{k,L}$ . If $k<0$, it can be shown that they span the whole space $H^+_{k,L}$ (see \cite{Br1} Proposition 1.12).
Moreover, one can check that $\xi_k(F_{h,m}(\tau,\kappa/2))$ equals up to a constant factor 
the usual holomorphic Poincar\'e series $P_{h,-m}(\tau)\in S_{\kappa,L^-}$ (as e.g.~defined in \cite{Br1} chapter 1.2.1). 
Therefore, for $\kappa >2$ (and $k <0$), the above duality statement also follows from chapter 1 of \cite{Br1}. However, the approach of the present paper is more conceptual while also covering the low weights $\kappa = 3/2$ and $2$ (where the $P_{h,-m}$ would need to be defined by Hecke summation).

The principal part of the Hejhal Poincar\'e series $F_{h,m}(\tau,\kappa/2)$ is equal to 
\[
e(m\tau)\frake_h +(-1)^q  e(m\tau)\frake_{-h}+\text{constant term}.
\]
This shows that, up to the constant term, any Fourier polynomial as in \eqref{principal} occurs as the principal part $P(f)$ of some $f\in H_{k,L}^+$. This can also be deduced from Theorem \ref{thmsurj} as follows:
\end{remark}

\begin{proposition}
For every Fourier polynomial of the form 
\[
Q(\tau)= \sum_{h\in L^\#/L}\sum_{\substack{n\in \Z+q(h) \\ n<0}} a^+(h,n) e(n\tau)\frake_h
\]
there exists an $f\in H_{k,L}^+$ with principal part $P(f)=Q+\frakc$ for some $T$-invariant constant  $\frakc\in \C[L^\#/L]$. The function $f$ is uniquely determined, if $k<0$ (and sometimes also if $k=0$).
\end{proposition}

\begin{proof}
Given $Q$ as above, we define a linear functional $\lambda_Q$ on $S_{\kappa,L^-}$ via the right hand side of (\ref{defpair'}).
By virtue of Corollary \ref{corpair} this functional is represented by a weak Maass form $h \in H_{k,L}^+$, i.e., $\{g,h\} = \lambda_Q(g)$ for all $g \in S_{\kappa,L^-}$.
The functional $\lambda_{Q-P(h)}$ then vanishes identically on $S_{\kappa,L^-}$.
Hence, by a variant of Theorem~3.1 of \cite{Bo2} (see \cite{Br1} Theorem~1.17), there exists a weakly holomorphic $h' \in M^!_{k,L}$ with principal part $Q-P(h)+\frakc$  for some $T$-invariant constant  $\frakc\in \C[L^\#/L]$.
Then $ f:=h+h' \in H^+_{k,L}$ has principal part $Q+\frakc$.
\end{proof}

\section{Special Schwartz Functions}

For the following discussion the reader should consult
\cite{KM90}, sections 5-8.

Kudla and Millson 
constructed (in more
generality) Schwartz forms $\varphi_{KM}$ on $V(\R)$ taking values
in $\mathcal{A}^q(D)$, the differential $q$-forms on $D$. More precisely,
\begin{equation*}
\varphi_{KM} \in [ S(V(\R)) \otimes \mathcal{A}^q(D)]^G \simeq
 [ S(V(\R)) \otimes {\bigwedge} ^q(\mathfrak{p^{\ast}})]^K,
\end{equation*}
where the isomorphism is again given by evaluation at the base point of
$D$.

We denote by $X_{\alpha\mu}$ ($ 1 \leq \alpha \leq p$, $p+1 \leq \mu \leq p+q$) the elements of the obvious basis of $\mathfrak{p}$ in (\ref{mathfrakp}) and let $\omega_{\alpha\mu}$ be the elements of the dual basis which pick out the $\alpha\mu$th coordinate of $\mathfrak{p}$.
Then $\varphi_{KM}$ is given by applying the operator
\begin{gather*}
\calD = \frac{1}{2^{q/2}} \prod_{\mu = p+1}^{p+q} \left[
\sum_{\alpha =1}^{p} \left(  x_{\alpha} - \frac{1}{2\pi}
\frac{\partial}{\partial x_{\alpha}}  \right)   \otimes
A_{\alpha\mu}  \right]
\end{gather*}
to the standard Gaussian $\varphi_0 \otimes 1 \in   [S(V(\R)) \otimes {\bigwedge} ^0(\mathfrak{p^{\ast}})]^K$:
\begin{equation*}
\varphi_{KM}= \calD (\varphi_0 \otimes 1).
\end{equation*}
Here $A_{\alpha\mu}$ denotes the left multiplication by $\omega_{\alpha\mu}$. Note that this is $2^{q/2}$ times the corresponding quantity in \cite{KM90}. It is easy to see that $\varphi_{KM}$ is $K$-invariant, and by \cite{KM1}, Theorem 3.1 it is an eigenfunction of $K'$ of weight $\tfrac{p+q}{2}$. 

For a multi-index $\underline{\alpha} =(\alpha_1,\dots, \alpha_q) \in \{1,\dots,p\}^q$, we have 
\[
\varphi_{KM}(x) = \sum_{\underline{\alpha}} P_{\underline{\alpha}}(x) \;\varphi_0(x) \; \omega_{\alpha_1 p+1} \wedge \dots \wedge \omega_{\alpha_q p+q}, 
\]
where $ P_{\underline{\alpha}}(x)$ is a (in general non-homogeneous) polynomial of degree 
$q$. For $\underline{\alpha} = (\alpha,\dots,\alpha)$, $ P_{\underline{\alpha}}(x)$ is given by 
\[
 P_{\underline{\alpha}}(x) = (4\pi)^{-q/2} H_q\left(\sqrt{2\pi}x_{\alpha}\right),
\]
where $H_q(t)= (-1)^q e^{t^2} \tfrac{d^q}{dt^q} e^{-t^2}$ is the $q$-th Hermite polynomial. For `mixed' $\underline{\alpha}$,  $P_{\underline{\alpha}}(x)$ is a product of Hermite functions in the $x_{\alpha}$.

In particular, we have (see \cite{KM1}, Prop.~5.1)
\begin{equation}\label{Euler}
\varphi_{KM}(0) = e_q,
\end{equation}
where for $q=2l$ even, $e_q$ is the Euler form of the symmetric space $D$ (which is the Euler class (see e.g.~\cite{KN}) of the tautological vector bundle over $D$, i.e., the fiber over a point $z \in D$ is given by the negative $q$-plane $z$) and zero for $q$ odd. Here $e_q$ is normalized such that it is given in $\bigwedge^q(\mathfrak{p}^{\ast})$ by
\begin{equation*}
e_q = \left(-\frac{1}{4\pi} \right)^l \frac{1}{l!} \sum_{\sigma \in S_q} \sgn(\sigma) \Omega_{p+\sigma(1),p+\sigma(2)} \dots  \Omega_{p+\sigma(2l-1),p+\sigma(2l)},
\end{equation*}
with
\begin{equation*}
\Omega_{\mu\nu} = \sum_{\alpha=1}^p \omega_{\alpha\mu} \wedge \omega_{\alpha\nu}.
\end{equation*}
Note that for $q=2$, $\Omega := - e_2$ is positive, i.e., it defines a K\"ahler form on the Hermitian domain $D$.

\medskip

  The space of $K'$-finite vectors in $\mathcal{S}(V(\R))$ is given by 
 the so-called polynomial Fock space $S(V(\R)) \subset \mathcal{S}(V(\R))$ which consists of those Schwartz functions on $V(\R)$ of the form $p(x)\varphi_0(x)$, where $p(x)$ is a polynomial function on $V(\R)$. 
Differentiating the action of $\Mp_2(\R) \times
\Orth(V(\R))$ on $S(V(\R))$ we obtain the associated action of the Lie
algebra $\mathfrak{sl}_2 \times \mathfrak{o}(V)$ which we also
denote by $\omega$. Then there is an intertwining map $\iota:
S(V(\R)) \longrightarrow \mathcal{P}(\C^{p+q})$ to the
infinitesimal Fock model of the Weil representation acting on the
space of complex polynomials $\mathcal{P}(\C^{p+q})$ in $p+q$
variables such that $\iota(\varphi_0) = 1$. We denote the
variables in $\mathcal{P}(\C^{p+q})$ by $z_{\alpha}$ ($1 \leq
\alpha \leq p$) and $z_{\mu}$ ($p+1 \leq \mu \leq p+q$). Since
\[
\iota
\left(  x_{j} - \frac{1}{2\pi} \frac{\partial}{\partial x_j}  \right)  \iota^{-1} =
 - \tfrac{1}{2\pi} z_{j},
\]
the form $\varphi_{KM}$ becomes
\begin{equation}\label{KM-Form}
\varphi_{KM} = \left(  \frac{-\sqrt{2}}{4\pi}   \right)^q \sum_{\alpha_1,\dots,\alpha_q} z_{\alpha_1} \cdots z_{\alpha_q} \otimes \omega_{\alpha_1 p+1} \wedge \cdots \wedge \omega_{\alpha_q p+q},
\end{equation}
when considered as an element in $[ \mathcal{P}(\C^{p+q}) \otimes \bigwedge^q(\mathfrak{p}^{\ast})]^K$.

In the Fock model, the Weil representation acts as follows: For $\mathfrak{o}(V(\R))$, we have
\begin{equation*}
\omega(X_{\alpha\mu}) = -4\pi \frac{\partial^2}{\partial z_{\alpha} \partial z_{\mu}} + \frac{1}{4\pi}z_{\alpha}z_{\mu}.
\end{equation*}
(There is a sign error in \cite{KM90}.) The elements $R$ and  $L$ in $\mathfrak{sl}_2(\C)$ (which correspond to the raising and lowering operators $R$ and $L$ on $\h$, see section 3) are given by
\begin{equation*}
R = \frac12
\begin{pmatrix}
1 & i \\ i & -1
\end{pmatrix}
\qquad \text{and} \qquad
L = \frac12
\begin{pmatrix}
1 & -i \\ -i & -1
\end{pmatrix},
\end{equation*}
and their action is
\begin{equation}\label{raising}
\omega(R) = - \frac{1}{8\pi} \sum_{\alpha=1}^p z^2_{\alpha} +2\pi \sum_{\mu = p+1}^{p+q} \frac{\partial^2}{\partial z_{\mu}^2},
\end{equation}
\begin{equation}\label{lowering}
\omega(L) =  2\pi \sum_{\alpha = 1}^{p} \frac{\partial^2}{\partial z_{\alpha}^2}-\frac{1}{8\pi} \sum_{\mu=p+1}^{p+q} z^2_{\mu}.
\end{equation}
The differentiation $d$ in the Lie algebra complex $[
\mathcal{P}(\C^{p+q}) \otimes \bigwedge(\mathfrak{p}^{\ast})]^K$
is given by
\begin{equation*}
d = \sum_{\alpha,\mu} \omega (X_{\alpha\mu}) \otimes
A_{\alpha\mu}.
\end{equation*}
Form this it is easy to see that $\varphi_{KM}$ is closed.


\medskip

We also define a form
\begin{equation}
\psi_{KM} = \psi \in [ \mathcal{P}(\C^{p+q}) \otimes {\bigwedge}
{^{q-1}}(\mathfrak{p}^{\ast})]^K
\end{equation}
(and therefore in $ [ S(V(\R)) \otimes \mathcal{A}^{q-1}(D) ]^G$)
by
\begin{equation}
\psi_{KM} = \frac{-1}{2(p+q-1)} h(\varphi_{KM}),
\end{equation}
where $h$ is the operator on $\left[ \mathcal{P}(\C^{p+q}) \otimes
\bigwedge ^{\ast}(\mathfrak{p}^{\ast})\right]^K$ given by
\begin{equation*}
h = \sum_{\alpha,\mu} z_{\mu} \frac{\partial}{\partial z_{\alpha}}
\otimes A^{\ast}_{\alpha\mu}.
\end{equation*}
Here $A^{\ast}_{\alpha\mu}$ denotes the (interior)
multiplication with $X_{\alpha\mu}$, i.e.,
$A^{\ast}_{\alpha\mu}(\omega_{\alpha'\mu'}) =
\delta_{\alpha\alpha'}\delta_{\mu\mu'}$. 

\begin{lemma}\label{psiformel}
\begin{equation*}
\psi_{KM} =  \left(  \frac{-\sqrt{2}}{4\pi}   \right)^q \frac{-1}{2 (q-1)!} \sum_{\alpha_1,\dots,\alpha_{q-1}} z_{\alpha_1} \cdots z_{\alpha_{q-1}} \det
\begin{pmatrix}
z_{p+1} &\dots & z_{p+q} \\
 \omega_{\alpha_1 p+1} & \dots & \omega_{\alpha_1 p+q} \\
 \vdots & & \vdots\\ 
\omega_{\alpha_{q-1} p+1} & \dots & \omega_{\alpha_{q-1} p+q}
\end{pmatrix},
\end{equation*}
where the determinant in the non-commutative algebra 
$\mathcal{P}(\C^{p+q}) \otimes \bigwedge^{\ast}(\mathfrak{p}^{\ast})$ is defined inductively via expansion after the first row.
\end{lemma}

\begin{proof}
We give a brief sketch. We first note (easily checked by induction)
\begin{equation}\label{phidet}
\varphi_{KM} =  \left(  \frac{-\sqrt{2}}{4\pi}   \right)^q \frac{1}{q!}\sum_{\alpha_1,\dots,\alpha_q} z_{\alpha_1} \cdots z_{\alpha_q} 
\det (\Omega),
\end{equation}
where
\[
\Omega=
\begin{pmatrix}
 \omega_{\alpha_1 p+1} & \dots & \omega_{\alpha_1 p+q} \\
 \vdots & & \vdots\\ 
\omega_{\alpha_{q} p+1} & \dots & \omega_{\alpha_{q} p+q}
\end{pmatrix}.
\]
The action of the operator $A^{\ast}_{\beta\mu}$ on $\det(\Omega)$ can be computed by expanding $\det (\Omega)$ by the $(p-\mu)$-th column. If we write $\Omega^{\nu\mu}$ for the $(q-1)\times (q-1)$ matrix obtained from $\Omega$ by canceling the $\nu$-th row and the $(p-\mu)$-th column, we have
\begin{align*}
A^{\ast}_{\beta\mu}\det(\Omega)&=-\sum_{\nu=1}^q(-1)^{\mu-p} A^{\ast}_{\beta\mu} \omega_{\alpha_\nu \mu} \det(\Omega^{\nu\mu})
=-\sum_{\nu=1}^q(-1)^{\mu-p} \delta_{\beta\alpha_\nu} \det(\Omega^{\nu\mu}) .
\end{align*}
A little calculation then shows
\begin{align*}
\psi_{KM}&=
  \frac{- (-\sqrt{2}/4\pi)^q}{2(p+q-1)(q-1)!}
\sum_{\alpha_1,\dots,\alpha_q} \frac{\partial}{\partial z_{\alpha_q}}
z_{\alpha_1} \cdots z_{\alpha_q} 
\det
\begin{pmatrix}
z_{p+1} &\dots & z_{p+q} \\
 \omega_{\alpha_1 p+1} & \dots & \omega_{\alpha_1 p+q} \\
 \vdots & & \vdots\\  
\omega_{\alpha_{q-1} p+1} & \dots & \omega_{\alpha_{q-1} p+q}
\end{pmatrix}.
\end{align*}
One now applies the product rule to obtain the assertion.
\end{proof}

\begin{theorem}[Kudla-Millson \cite{KM90}]\label{KMpsi}\quad
\begin{itemize}
\item[(i)] We have the identity
\[
\omega(L)   \varphi_{KM} =  d \psi_{KM}.
\]

\smallskip

\item[(ii)]

$\psi_{KM}$ is an eigenfunction of $K'$ with weight $\tfrac{p+q}2
-2$.
\end{itemize}
\end{theorem}

\begin{example}
For signature $(p,1)$, we have in the Fock model (Schr\"odinger model) of the Weil
representation
\[
\varphi_{KM} = -\frac{\sqrt{2}}{4\pi} \sum_{\alpha=1}^{p}
z_{\alpha} \otimes  \omega_{\alpha p+1} = (\iota \otimes 1) \left( \sqrt{2}
\sum_{\alpha=1}^p x_{\alpha} \varphi_0 \otimes \omega_{\alpha p+1}
\right)
\]
and
\[
\psi_{KM} = \frac{\sqrt{2}}{8 \pi} z_{p+1} \otimes 1 = (\iota \otimes 1) \left(
\frac{-1}{\sqrt{2}} x_{p+1} \varphi_0  \otimes 1 \right).
\]
\end{example}

\bigskip

We now specialize to the case of signature $(p,2)$, i.e., $q=2$,
when we have an underlying complex structure on $D$. Then by the
conventions of section 2, the operators $\partial$ and
$\bar{\partial}$ on $D$ are given by
\begin{gather*}
\partial = \frac12 \sum_{\alpha} \omega(X_{\alpha p+1} - i X_{\alpha p+2}) \otimes ( \omega_{\alpha p+1} + i\omega_{\alpha p+2}), \\
\bar{\partial} = \frac12 \sum_{\alpha} \omega(X_{\alpha p+1} + i X_{\alpha p+2}) \otimes ( \omega_{\alpha p+1} - i\omega_{\alpha p+2}).
\end{gather*}

We set $d^c = \tfrac{1}{4\pi i}(\partial - \bar{\partial})$, so that $dd^c = - \tfrac{1}{2\pi i} \partial \bar{\partial}$. Note that in this case $\varphi_{KM}$ is actually a $(1,1)$-form.

\begin{theorem}\label{ddcGleichung}
In the case of signature $(p,2)$, we have
\[
\omega(L) \varphi_{KM} = -dd^c \varphi_0.
\]
\end{theorem}

We prove the theorem by computing the two sides separately. For
the left hand side, by (\ref{lowering}) and (\ref{KM-Form}) we
easily see
\begin{lemma} For signature $(p,2)$,
\[
\omega(L) \varphi_{KM} =  \frac{1}{2\pi}  \sum_{\alpha=1}^p
\omega_{\alpha p+1} \wedge \omega_{\alpha p+2} - \frac1{64 \pi^3}
(z^2_{p+1} + z^2_{p+2}) \sum_{\alpha_1,\alpha_2}   z_{\alpha_1}
z_{\alpha_2} \, \omega_{\alpha_1 p+1} \wedge \omega_{\alpha_2
p+2}.
\]
\end{lemma}

For the right hand side we get the same:
\begin{lemma} For signature $(p,2)$,
\[
-dd^c \varphi_0 =  \frac{1}{2\pi}  \sum_{\alpha=1}^p \omega_{\alpha
p+1} \wedge \omega_{\alpha p+2} - \frac1{64 \pi^3} (z^2_{p+1} +
z^2_{p+2}) \sum_{\alpha_1,\alpha_2} z_{\alpha_1} z_{\alpha_2}  \,
\omega_{\alpha_1 p+1} \wedge \omega_{\alpha_2 p+2}.
\]
\end{lemma}

\begin{proof}
We have
\begin{equation*}
 \frac12 \omega(X_{\alpha p+1} \pm i X_{\alpha p+2}) = -2\pi \frac{\partial}{\partial z_{\alpha}}
( \frac{\partial}{\partial z_{p+1}} \pm i \frac{\partial} {\partial z_{p+2}}) +\frac{1}{8\pi}z_{\alpha}(z_{p+1} \pm i z_{p+2}).
\end{equation*}
We first see
\begin{equation*}
\bar{\partial}(\varphi_0) = \frac{1}{8\pi} \sum_{\alpha} z_{\alpha}(z_{p+1} + iz_{p+2}) (\omega_{\alpha p+1} -i \omega_{\alpha p+2})
\end{equation*}
 and then
\begin{align*}
\partial  \bar{\partial}(\varphi_0) &= -\frac12 \sum_{\alpha}
 (\omega_{\alpha p+1} +i \omega_{\alpha p+2}) \wedge (\omega_{\alpha p+1} +i \omega_{\alpha p+2})  \\
&\quad - \frac{1}{64\pi^2} \sum_{\alpha_1,\alpha_2} z_{\alpha_1} z_{\alpha_2}
 (z_{p+1} - iz_{p+2})(z_{p+1} + iz_{p+2}) \notag \\
& \qquad \qquad \qquad \qquad \times
(\omega_{\alpha_1 p+1} +i \omega_{\alpha_1 p+2}) \wedge (\omega_{\alpha_2 p+1} -i \omega_{\alpha_2 p+2}) \notag \\
&=
 +i  \sum_{\alpha}
 \omega_{\alpha p+1} \wedge \omega_{\alpha p+2}  \\
&\quad -  \frac{i}{32\pi^2} \sum_{\alpha_1,\alpha_2} z_{\alpha_1} z_{\alpha_2}
(z_{p+1}^2 + z^2_{p+2})  \, \omega_{\alpha_1 p+1} \wedge \omega_{\alpha_2 p+2}. \notag
\end{align*}
The lemma follows.
\end{proof}

\begin{remark}
 Theorem \ref{KMpsi} states that there is $(q-1)$-form $\psi$ such that
$\omega(L) \varphi_{KM} = d \psi_{KM}$ (for any signature). So Theorem
\ref{ddcGleichung} suggests (and it is easily checked) that in the Hermitian case one has
\begin{equation}\label{dcformel}
\psi = - d^c \varphi_0.
\end{equation}
From this perspective we can consider
Theorem \ref{ddcGleichung} as the Weil representation-theoretic
analogue of the $dd^c$-Lemma in complex geometry.
\end{remark}

\section{The Theta Lifts}

\subsection{The Kudla-Millson Lift}
\label{sect6.1}
We consider the Schwartz function $\varphi_{KM}$ from the previous section. Since $\varphi_{KM}$ is closed and an eigenfunction of $K'$ of weight $\kappa=(p+q)/2$, we get
\begin{equation*}
\Theta(\tau,z,\varphi_{KM}) \in {A}_{\kappa,L} \otimes
\mathcal{Z}^q(X),
\end{equation*}
i.e., $\Theta(\tau,z,\varphi_{KM})$ is a non-holomorphic modular
form of weight $\kappa$ associated to the representation $\varrho_L$
with values in $\mathcal{Z}^q(X)$, the closed $q$-forms on $X$. (In fact, it extends to a closed $q$-form on the Borel-Serre compactification $\bar{X}$ of $X$, see \cite{FM1}).
We then can consider the assignment
\begin{equation*}
\eta \mapsto \Lambda_{KM}(\tau,\eta) := \int_X \eta \wedge
\Theta(\tau,z,\varphi_{KM}),
\end{equation*}
which is defined for rapidly decreasing $(p-1)q$-forms on $X$.
This map factors through $H_c^{(p-1)q}(X)$, the de Rham
cohomology with compact support of $X$.

\begin{theorem}[Kudla-Millson, see \cite{KM90}]\label{KMTheorem}
For $\eta$ closed and rapidly decreasing, the function $\Lambda(\tau,\eta)$ is 
\emph{holomorphic} on $\H$ and at the cusp (even though $
\Theta(\tau,z,\varphi_{KM})$ is not). We therefore have a map
\begin{equation*}
\Lambda_{KM}: H_c^{(p-1)q}(X) \longrightarrow M_{\kappa,L}.
\end{equation*}
The Fourier expansion is given by
\begin{equation*}
\Lambda_{KM}(\tau,\eta)_h = \delta_{h0} \, \left(\int_X \eta \wedge e_q \right) +
\sum_{n> 0} \left( \int_{Z(h,n)} \eta \right) \, e(n \tau),
\end{equation*}
i.e., the $n$-th Fourier coefficient of
$\Theta(\tau,z,\varphi_{KM})_h$ is a Poincar\'e dual form of the
cycle $Z(h,n)$. Here (for $q$ even), $e_q$ is the Euler form of the symmetric space $D$, see (\ref{Euler}). For $q$ odd, $\Lambda_{KM}(\tau,\eta)$ is a cusp form.
\end{theorem}

\begin{remark}

For $p=1$ and $X$ non-compact, we have $H_c^0(X)=0$, so that Theorem \ref{KMTheorem} would be empty. However, in that case, the $q$-form $\Theta(\tau,z,\varphi_{KM})$ is actually rapidly decreasing on $X$, see \cite{FM1}, and therefore defines a class in $H^q_c(X)$. This corresponds to the fact that the cycles $Z(h,n)$ are a collection of points in this situation.
Hence $\Lambda_{KM}$ is defined on $H^0(X) \simeq \R$ and for $q \geq 3$, the lift $\Lambda(\tau,\eta)$ is a holomorphic modular form (\cite{FM1}), while for $q=2$, it is \emph{non}-holomorphic, see \cite{Fu}.

\end{remark}

\subsection{The (Generalized) Borcherds Lift}

\subsubsection{The Hermitian Case}
Here we assume that $V$ has signature $(p,2)$ such that the corresponding symmetric domain is Hermitean.
We put $k=1-p/2$ and $\kappa=1+p/2$.  
In this section we briefly recall some facts on the Borcherds lift from weakly holomorphic elliptic modular forms of weight $k$ to automorphic forms on the orthogonal group $\Orth(V(\R))$ \cite{Bo1} and its generalization to weak Maass forms \cite{Br1}.


We consider the Siegel theta function
$\Theta(\tau,z,\varphi_0^{p,2}) \in {A}_{-k,L} \otimes
C^{\infty}(X)$.
The (additive) Borcherds lift of a weak Maass form $f\in H^+_{k,L^-}$ is defined by the theta integral
\begin{equation}\label{defbor}
\Phi(z,f) = \Phi(z,f,\varphi_0) = \int_{\G' \back \h}^{reg} \langle f(\tau), \overline{\Theta(\tau,z,\varphi_0^{p,2})} \rangle \,d \mu.
\end{equation}
The integral is typically divergent and needs to be regularized in the following way (indicated by the superscript ``reg''): If $F$ is a $\G'$-invariant function on $\h$, we consider for an additional complex parameter $s\in \C$ the function
\begin{equation}\label{petreg}
\lim_{t\to\infty} \int_{\calF_t} F(\tau) v^{-s} \,d\mu,
\end{equation}
where $\calF_t$ denotes the truncated fundamental domain (\ref{Fdomain}) for the action of $\Sl_2(\Z)$ on $\H$. Formally, at $s=0$, the quantity \eqref{petreg} equals the usual integral of $F$ over $\G' \back \h$.
However, even if this diverges, \eqref{petreg} sometimes converges for $\Re(s)\gg 0$ and has a meromorphic continuation to a neighborhood of $s=0$.
Then we define 
\begin{equation}\label{petreg2}
\int^{reg}_{\G' \back \h} F(\tau) \,d\mu
\end{equation}
to be the constant term in the Laurent expansion of \eqref{petreg} at $s=0$.
So one feature of the regularization  consists in prescribing the order of integration, the other in introducing $s$ and looking at the Laurent expansion.
The above regularization of the theta integral for weakly holomorphic modular forms was discovered  by Harvey and Moore \cite{HM} and vastly generalized by Borcherds \cite{Bo1}. In connection with weak Maass forms  it was investigated in \cite{Br1}.

\begin{remark}
Note that in \cite{Bo1} and \cite{Br1} the Borcherds lift is actually defined for signature $(2,p)$. Identifying the symmetric spaces $D_{p,2}$ and $D_{2,p}$ and switching to the space $V^-$ one has
\[
\Phi(z,f) = \left( f(\tau), \, \Theta(\tau,z,\varphi_0^{2,p},L^-)\right)^{reg}_{k,L^-},
\]
where the integral defining the Peterson scalar product is regularized as above. The right hand side is the definition used in \cite{Bo1} and \cite{Br1}. The definition given here is more convenient for our purposes.
\end{remark}

\begin{proposition}\label{propreg}
Denote the Fourier expansion of $f\in H^+_{k,L^-}$ as in \eqref{deff}.
The regularized integral for $\Phi (z,f)$ converges to a $\G$-invariant $C^\infty$ function on $D$ with a logarithmic singularity along the divisor
\[
-2\sum_h\sum_{n\in\Q_{<0}} a^+(h,n) Z(h,-n).
\]  
\end{proposition}

\begin{proof}
See \cite{Bo1} section 6 and  \cite{Br1} section 2.2. 
\end{proof}


We list some further properties of $\Phi (z,f)$:
Outside its singularities, the function $\Phi (z,f)$ is almost an eigenfunction of the invariant Laplacian $\Delta$ on $D$ induced by the Casimir element of the universal enveloping algebra of the Lie algebra of $\Orth(V(\R))$.
More precisely, if we normalize $\Delta$ as in \cite{Br1} chapter 4.1, we have $\Delta \Phi (z,f) = \frac{p}{4} a^+(0,0)$.
Consequently, since $\Delta$ is an elliptic differential operator on $D$, the function $\Phi (z,f)$ is actually real analytic outside its singularities. 

It is proved in \cite{Br1} chapter 3.2 that $\Phi (z,f)$ can be split 
into a sum
\begin{equation}\label{splitequ}
\Phi (z,f)=-2 \log|\Psi(z,f)| + \xi(z,f),
\end{equation}
where $\xi(z,f)$ is real analytic on the whole domain $D$ and $\Psi(z,f)$ is a meromorphic function on $D$ whose divisor equals 
\begin{equation}\label{eqdiv}
Z(f):=\sum_h\sum_{n\in\Q_{<0}} a^+(h,n) Z(h,-n).
\end{equation}
If $f \in M^{!}_{k,L^-}$ and its Fourier coefficients $a^+(h,n)$ with negative index are integral, then $\xi(z,f)$ reduces to a ``simple'' function and  $\Psi(z,f)$ is a meromorphic modular form for the group $\Gamma$ of weight $a^+(0,0)$. This is the (multiplicative) Borcherds lift from $M_{k,L^-}^!$ to meromorphic modular forms for $\Gamma$. If $f$ is not holomorphic, then $\Psi(z,f)$ is far from being modular under $\Gamma$.
(Caution: It is in general not true that $\xi(z,f)$ is equal to the regularized theta integral of $f^-$ as one might think.) 


From these considerations it can be deduced that the $(1,1)$-form
\begin{equation}\label{deflb}
\Lambda_B(f):= dd^c \Phi (z,f)= dd^c \xi(z,f)
\end{equation}
is closed, harmonic,  $\Gamma$-invariant, and of moderate growth. It represents the Chern class of the divisor \eqref{eqdiv} 
in the second cohomology $H^2(X)$.  
In particular, we have 
\begin{equation}\label{pdbor}
\int_{Z(f)} \eta = \int_X \eta\wedge \Lambda_B(f)
\end{equation}
for any rapidly decreasing $2(p-1)$-form $\eta$. 
If $f \in M^!_{k,L^-}$, then $\Lambda_B(f) = a^+(0,0) \Omega$, i.e., a multiple of the K\"ahler form $\Omega= - e_2$ on $D$ and actually vanishes for $f \in M^{!!}_{k,L^-}$ and $k<0$.

We summarize part of the above discussion in the following theorem.

\begin{theorem}\label{thmlb} (See \cite{Br1} chapter 5.)
The assignment $f\mapsto \Lambda_B(f)$ defines a linear mapping $H^+_{k,L^-}\to\calZ^{1,1}_{mg}(X)$ to the space $\calZ^{1,1}_{mg}(X)$ of closed $(1,1)$-forms of moderate growth on $X$. The induced mapping $H^+_{k,L^-}\to\calZ^{1,1}_{mg}(X)/\C\Omega$ factors through $H^+_{k,L^-}/ M^!_{k,L^-}$. 
\end{theorem} 


\subsubsection{The Real Case}

We now consider the case of general signature $(p,q)$. If $x\in V(\R)$, we put $|x|=|(x,x)|^{1/2}$. Moreover, as before we let $\kappa = (p+q)/2$ and $k = 2-\kappa$.
The  Schwartz form  $\psi= \psi_{KM}$ (see section 4) has
weight $-k= (q+p)/2 -2 $ by Theorem \ref{KMpsi} (ii), hence
\begin{equation*}
\Theta(\tau,z,\psi) \in A_{-k,L} \otimes
\mathcal{A}^{q-1}(X).
\end{equation*}
Note that the components of
$\Theta(\tau,z,\psi)$ are given by
\begin{equation*}
\left(\Theta(\tau,z,\psi)\right)_h = v^{1-q/2}
\sum_{\la \in L+h} P_{\psi}(\sqrt{v}\la,z) \varphi^{p,q}_0(\la,\tau,z),
\end{equation*}
with a (in general non-homogeneous) polynomial $P_{\psi}(x,z) \in
[\mathcal{P}(V(\R)) \otimes \mathcal{A}^{q-1}(D)]^G$ of degree
$q$. By Lemma~\ref{psiformel}, $P_{\psi}$ can be written as
\begin{equation}
P_{\psi}(x,z) = \sum_{\underline{\alpha}} Q_{\underline{\alpha}}(x,z) \, R_{\underline{\alpha}}(x,z),
\end{equation}
where the sum extends over all multi-indices $\underline{\alpha} = (\alpha_1,\dots, \alpha_{q-1})$. Here $Q_{\underline{\alpha}}$ is polynomial in $[\mathcal{P}(V(\R)) \otimes \mathcal{A}^{0}(D)]^G$ of degree $q-1$ (in fact, a product of Hermite functions as for $\varphi_{KM}$), which only depends on $x_{z^{\perp}}$; while $
 R_{\underline{\alpha}}$ is a homogeneous linear polynomial in $[\mathcal{P}(V(\R)) \otimes \mathcal{A}^{q-1}(D)]^G$, which only depends on $x_z$.
In particular, we have $P_{\psi}(0,z)=0$. We normalize the polynomials such that at the base point $z_0$, we have
\[
R_{\underline{\alpha}}(x) =
R_{\underline{\alpha}}(x,z_0) 
= \frac{-1}{\sqrt{2}(q-1)!} \det 
\begin{pmatrix}
       x_{p+1} & \cdots & x_{p+q} \\
       \omega_{\alpha_1 p+1} &\cdots & \omega_{\alpha_{1} p+q} \\
       \vdots & & \vdots\\
       \omega_{\alpha_{q-1} p+1} &\cdots & \omega_{\alpha_{q-1} p+q }
\end{pmatrix},
\]
and the leading term of  $Q_{\underline{\alpha}}(x) = Q_{\underline{\alpha}}(x,z_0)$ is given by $2^{(q-1)/2} \prod_{i=1}^{q-1} x_{\alpha_i}$.

Let $f \in H^+_{k,L^-}$ be a weak Maass form. We then define the Borcherds lift for signature $(p,q)$ by
\begin{equation}\label{genborch}
\Phi(z,f,\psi) := \left(\Theta(\tau,z,\psi), v^k\overline{f}\right)^{reg}_{-k,L} = 
\int_{\G' \back \h}^{reg} \langle {\Theta(\tau,z,\psi)}, \overline{f} \rangle \,d\mu.
\end{equation}


We define the cycle $Z(f)$ associated to $f$ by (\ref{eqdiv}), as in the Hermitian case. Note that the Fourier coefficients satisfy $a^+(-h,n) = (-1)^q a^+(h,n)$, corresponding to $Z(-h,n) = (-1)^q Z(h,n)$. 
In particular, we have
$ Z(F_{h,m}) = 2 Z(h,-m) $,
where $F_{h,m}$ is the Hejhal Poincar\'e series, see Remark~\ref{Hejhal}.

\begin{proposition}\label{sing}
The regularized integral for $\Phi(z,f,\psi)$ converges to a $\G$-invariant $C^\infty$ function on $D$ with singularities along $Z(f)$. More presicely, there exist $G$-invariant scalar-valued polynomials $\widetilde{Q}_{\underline{\alpha}}(x,z) \in [\mathcal{P}(V(\R)) \otimes \mathcal{A}^{0}(D)]^G$ of degree $(q-1)$ such that in a small neighborhood 
of $ w \in D$ the singularity is of type
\[
\sum_{h}\sum_{n \in \Q_{<0}} a^+(h,n)
\sum_{\lambda \in w^\perp \cap L_{h,-n} }  
\sum_{\underline{\alpha}} 
\widetilde{Q}_{\underline{\alpha}}\left( \lambda_{z^{\perp}}/ |\lambda_z|,z 
\right)
 R_{\underline{\alpha}}\left( \lambda_z/|\lambda_z|,z\right). 
\]
The leading term   $\widetilde{Q}_{\underline{\alpha},q-1}$ of  $\widetilde{Q}_{\underline{\alpha}}$ at the base point $z_0$ is given by
\[
\widetilde{Q}_{\underline{\alpha},q-1}(x) = 
 \frac{\G\left(q/2\right)}{\sqrt{2}\pi^{q/2}} \prod_{i=1}^{q-1} x_{\alpha_i}.
\]
\end{proposition}

Note that all sums in the above formula are finite. The polynomials $\widetilde{Q}_{\underline{\alpha}}$ and $R_{\underline{\alpha}}$ only depend on the signature $(p,q)$ and not on $f$.

\begin{example}\label{hyper}
In the hyperbolic case ($q=1$), we have  $\widetilde{Q} =1/\sqrt{2}$ and $R(x,z)= \frac{1}{\sqrt{2}}(x,gv_{p+q}) $ with $z=gz_0$ for some $g \in G$, and  
the singularity near $w\in D$ is given by
\[
\frac{1}{2} 
\sum_{h}\sum_{n \in \Q_{<0}} a^+(h,n)
\sum_{\lambda \in w^\perp \cap L_{h,-n} }  
\sum_{\underline{\alpha}} \sgn\left(\la,g v_{p+q}\right). 
\]
In particular, $\Phi(z,f,\psi)$ is locally bounded.
\end{example}

\begin{proof}
The argument follows in large parts \cite{Bo1}, section 6.
Since the function $f^-$ is exponentially decreasing as $v\to \infty$, and 
$\Theta(\tau,z,\psi)$ has moderate growth, the  integral 
\[
\int_{\calF} \langle f^-(\tau), \overline{\Theta(\tau,z,\psi)} \rangle \,
d\mu
\]
over the standard fundamental domain $\calF= \{ \tau\in \H;\; \text{$|\tau|\geq 1$, $|u|\leq 1/2$}\}$ for $\Sl_2(\Z)$ converges absolutely and defines a real analytic function on $D$.
Moreover, the integral of $\langle f^+(\tau), \overline{\Theta(\tau,z,\psi)} \rangle $ over the compact subset $\calF_1\subset\H$ converges absolutely and defines a real analytic function on $D$. Hence we are left with considering the function
\[
h(z,s) =\int_1^\infty \int_{-1/2}^{1/2} \langle f^+(\tau), \overline{\Theta(\tau,z,\psi)} \rangle v^{-2-s}\,du\,dv.
\] 
If we insert the Fourier expansions, the integration over $u$ picks out the constant term in the Fourier expansion of $\langle f^+(\tau), \overline{\Theta(\tau,z,\psi)} \rangle$. Thus
\begin{align}\label{c1}
h(z,s) &=\sum_{\lambda\in L^\#}
a^+(\lambda,-q(\lambda))  \int_1^\infty P_{\psi}(\sqrt{v}\la,z)\exp( 4\pi v q(\lambda_{z})) v^{-1-s}dv.
\end{align}
(Note that  $q(\lambda_{z})\leq 0$ for all $\lambda$ in the above sum.)
By the growth of  $f$, it has only finitely many non-zero Fourier coefficients $a^+(h,n)$ with negative index $n$.
Hence, if $U\subset D$ is a relatively compact neighborhood of $w$ and $\eps>0$, then by reduction theory the set
\[
S_f(U,\eps)=\left\{ \lambda\in L^\#;\quad 
\text{$a^+(\lambda,-q(\lambda))\neq 0$ and $|q(\lambda_{z})|<\eps$ for some $z\in U$}\right\}
\]
is finite. We split the sum over $\lambda\in L^\#$ in \eqref{c1} into the sum over $\lambda\in S_f(U,\eps)$ and the sum over $\lambda\in L^\#-S_f(U,\eps)$.
The latter sum is up to a constant factor majorized by
\[
\sum_{\lambda\in L^\#-S_f(U,\eps)} 
|a^+(\lambda,-q(\lambda))|  \exp( 2\pi  q(\lambda_{z})),
\]
locally uniformly in $s\in\C$ and $z\in U$.
According to Lemma \ref{growth2} there exists a constant 
$C>0$ such that this sum is majorized by
\begin{align*}
&\sum_{\substack{\lambda\in L^\#\\ a^+(\lambda,-q(\lambda))\neq 0}}
\exp\left( C\sqrt{|q(\lambda)|}+2\pi  q(\lambda_{z})\right) \\
&= \sum_{\substack{\lambda\in L^\#\\ a^+(\lambda,-q(\lambda))\neq 0}} 
\exp\left( C\sqrt{|q(\lambda)|}+\pi q(\lambda)\right)
\exp\left(-\frac{\pi}{2} (\la,\la)_z\right)
\ll \sum_{\substack{\lambda\in L^\#}} \exp\left(-\frac{\pi}{2} (\la,\la)_z\right),
\end{align*}
which clearly converges.
We may conclude that the sum over $\lambda\in L^\# - S_f(U,\eps)$ in \eqref{c1} (and all its derivatives) converges locally uniformly absolutely for $s\in \C$ and $z\in U$. In particular, it defines a $C^\infty$ function for $z\in U$ at $s=0$.

Hence the singularity of $\Phi (z,f,\psi)$ for $z\in U$ is dictated by the finite sum
\begin{equation*}
h_1(z,s) =\sum_{\lambda\in S_f(U,\eps)}
a^+(\lambda,-q(\lambda))  \int_1^\infty P_{\psi}(\sqrt{v}\la,z) \exp( 4\pi v q(\lambda_{z})) v^{-1-s}dv.
\end{equation*}
Note that so far the argument works for any polynomial $P(x,z)$.

Since $P_{\psi}(0,z) =0$, the term for $\la =0$ does not contribute (and this implies that $h_1(z,s)$ will be actually holomorphic at $s=0$). We are left to consider (for $\la \ne 0$)
\begin{equation}\label{sing1}
\begin{split}
\tilde{\psi}(\lambda,z,s)& :=\int_1^\infty P_{\psi}(\sqrt{v}\la,z) \exp( 4\pi v q(\lambda_{z})) 
v^{-1-s}dv  \\
& = \sum_{\underline{\alpha}} R_{\underline{\alpha}}(\la_z,z) 
\int_1^\infty Q_{\underline{\alpha}}(\sqrt{v}\la_{z^{\perp}},z) \exp( 4\pi v q(\lambda_{z})) v^{1/2-s} \frac{dv}{v}  \\
& = \sum_{\underline{\alpha}} R_{\underline{\alpha}}(\la_z,z) \sum_{\ell=0}^{q-1}  
Q_{\underline{\alpha},\ell}(\la_{z^{\perp}},z)
\int_1^\infty \exp( 4\pi v q(\lambda_{z})) v^{\ell/2+1/2-s} \frac{dv}{v}  \\
& = \sum_{\underline{\alpha}} R_{\underline{\alpha}}(\la_z,z) \sum_{\ell=0}^{q-1}  
Q_{\underline{\alpha},\ell}(\la_{z^{\perp}},z) \left|4\pi q(\la_z) \right|^{s-(\ell+1)/2} \G\left(\tfrac{\ell+1}{2}-s,|4\pi q(\la_z)|\right). 
\end{split}
\end{equation}
Here $Q_{\underline{\alpha},\ell}(\la_{z^{\perp}},z)$ denotes the homogeneous component of degree $\ell$ of $Q_{\underline{\alpha}}(\la_{z^{\perp}},z)$ and 
 $\Gamma(a,x)=\int_x^\infty e^{-t} t^{a-1} dt$ is the incomplete Gamma function as in \cite{AS} (6.5.3). This shows that $h_1(z,s)$ has a meromorphic continuation in $s$ to the whole complex plane and is holomorphic at $s=0$.
The claim now follows from the recurrence relations of the incomplete Gamma function and from $\G(q/2,|4\pi q(\la_z)|) = \G(q/2) + O( |\la_z|)$ as $\la_z \to 0$.
%
%
\end{proof}

Analogously to the Hermitian situation we define the map $\La_{B,\psi}$ by
\begin{equation}
\Lambda_{B,\psi}(z,f) = - d \Phi(z,f,\psi).
\end{equation}
This is the analogue to the map $\Lambda_B = \Lambda_{B,\varphi_0}$ considered in the Hermitian case. The corresponding geometric properties will easily follow from the relationship to the Kudla-Millson lift which we establish in the next section.

\begin{remark}
Note that in the Hermitian case of signature $(p,2)$ we have now defined two Borcherds lifts, the `classical' one $\Phi(z,f) = \Phi(z,f,\varphi_0)$ and 
$\Phi(z,f,\psi)= \Phi(z,f,- d^c\varphi_0)$. We will later see that we actually have (outside $Z(f)$)
\[
d^c \Phi(z,f) =  - \Phi(z,f,\psi);
\]
i.e., one can interchange the differential operator $d^c$ on the orthogonal domain $D$ with the regularized $\Sl_2(\R)$-integral.
\end{remark}

\begin{remark}\label{harmonic}
Analogously to \cite{Br1}, chapter 4, one can show that $\Phi(z,f,\psi)$ is
harmonic, i.e., annihilated by the Laplacian for the symmetric space $D$. In particular, $\Phi(z,f,\psi)$ is real analytic.
\end{remark}

\section{Main Results}
\label{sect7}

Recall (with  $\kappa = (p+q)/2$ and $k= 2- \kappa$) the map 
$\xi_k(f) = R_{-k}(v^k \bar{f}) \in S_{\kappa,L}$ for
$f \in H^+_{k,L^-}$.
Since $\Theta(\tau,z,\varphi_{KM})$ is moderately increasing in $\tau$, we can therefore study the scalar product 
\begin{equation*}
(\Theta(\tau,z,\varphi_{KM}), \xi_{k}(f))_{\kappa,L} \in \mathcal{Z}^{q}(X).
\end{equation*}

\begin{theorem} \label{Main1}
Let $V$ be of signature $(p,q)$ and let $f \in H^+_{k,L^-}$. Then, outside $Z(f)$, the set of singularities of the Borcherds lift $\Phi(z,f,\psi)$, we have
\begin{equation}\label{Main1a}
\left( \Theta(z,\varphi_{KM}), \xi_k(f)\right )_{\kappa,L} -
 a^+(0,0) e_q  = - d\Phi(z,f,\psi) =  \Lambda_{B,\psi}(z,f).
\end{equation}
For the `classical' Borcherds lift in signature $(p,2)$, we have, outside $Z(f)$,
\begin{equation}\label{Main1b}
\left(\Theta(z,\varphi_{KM}), \xi_k(f) \right)_{\kappa,L}  + a^+(0,0) \Omega = dd^c 
\Phi(z,f) =  \Lambda_B(z,f).
\end{equation}
Here $a^+(0,0)$ is the constant coefficient of the Fourier expansion of $f$, $e_q$ is the Euler form on $D$ for $q$ even (and zero for $q$ odd), and $\Omega = - e_2$ is the K\"ahler form for the Hermitian domain $D_{p,2}$. 

\end{theorem}

\begin{theorem}\label{BorcherdsProp}\quad

\begin{itemize}
\item[(i)]
$\Lambda_{B,\psi}(f)$ extends to a smooth closed $q$-form of moderate growth on $X$.

\item[(ii)]
If  $f \in M^!_{k,L^-}$ is weakly holomorphic, then
\[
\Lambda_{B,\psi}(f) = - a^+(0,0) e_q.
\]

\item[(iii)]
 In particular, $\Lambda_{B,\psi}$ induces a map
\[
\Lambda_{B,\psi}: H^+_{k,L^-} / M^!_{k,L^-} \longrightarrow H^q(X)/ \C e_q.
\]

\end{itemize}
\end{theorem}

\begin{proof}
This follows immediately from Theorem~\ref{Main1}: The left hand side of (\ref{Main1a}) is smooth and has moderate growth, which gives (i). The assertion (ii) follows from $\xi_k(f) =0$ for $f \in M^!_{k,L^-}$.
\end{proof}

 Observe that in view of (\ref{Main1b}) one obtains new proofs of the analogous properties of the `classical' Borcherds lift listed in section 6.

\medskip

We denote by $(\cdot,\cdot)_X$ the natural pairing between $\calZ_{rd}^{qp-\ell}(X)$, the closed rapidly decreasing $(qp-\ell)$-forms, and $\calZ_{mg}^{\ell}(X)$, the closed $\ell$-forms of moderate growth. Recall that
by Poincar\'e duality this induces a non-degenerate pairing between ${H}_c^{qp-\ell}(X)$ and ${H}^{\ell}(X)$, see \cite{BT}.

On the other hand, recall the pairings $\{\,,\,\}$ and $\{\,,\,\}'$ between $M_{\kappa,L}$ and $H^+_{k,L^-}$ defined in \eqref{defpair} and \eqref{defpair'}.


\begin{theorem}\label{Main2}
The Kudla-Millson lift $
\Lambda_{KM}: \calZ_{rd}^{q(p-1)}(X) \to M_{\kappa,L}$ and the Borcherds lift
$\Lambda_{B,\psi}: H^+_{k,L^-} \to \calZ^q_{mg}(X)$ are adjoint in the following sense: For any $\eta\in \calZ_{rd}^{q(p-1)}(X)$ and any $f\in H^+_{k,L^-}$ we have
\[
\big(\eta, \Lambda_{B,\psi}(f)\big)_X  =\left\{\Lambda_{KM}(\eta),f\right\}'.
\]
Moreover, this duality factors through the non-degenerate pairings on the cohomology level and the one of $M_{\kappa,L}$ with  $H^+_{k,L^-} / M^{!!}_{k,L^-}$.
\end{theorem}

\begin{proof} 
Using \eqref{Main1a} we see that
\begin{align*}
 \big(\eta, \Lambda_{B,\psi}(f)\big)_X +  a^+(0,0) \big(\eta, e_q\big)_X
&=\big(\eta, \,(\Theta_{KM}(\tau,Z),\xi_k(f))_{\kappa,L} \big)_X\\
&=\big( (\eta,\Theta_{KM}(\tau,Z))_X, \,\xi_k(f) \big)_{\kappa,L}\\
& = \big\{\Lambda_{KM}(\eta),f\big\}.
\end{align*}
Here, in the second equality, we have exchanged the order of integration. This is valid because the latter integrals converge absolutely. Now the assertion follows from the Fourier expansion of the Kudla-Millson lift (see Theorem \ref{KMTheorem}) and the definition of the pairing $\{\,,\,\}'$.
\end{proof}

\begin{remark}
In signature $(p,2)$, the same statement of Theorem~\ref{Main2} holds with $\Lambda_{B,\psi}$ replaced by the classical Borcherds lift $\Lambda_B$.
\end{remark}

\begin{corollary}\label{c99}
For $\eta \in \calZ_{rd}^{q(p-1)}(X)$ and $f\in H_{k,L^-}^+$, we have
\[
\big(\eta, \Lambda_{B,\psi}(f)\big)_X = \int_{Z(f)} \eta.
\]
For $q$ even, this also holds with $\eta = e_q^{p-1}$. In particular, 
$\Lambda_{B,\psi}(f)$ is a  harmonic representative of the Poincar\'e dual form of the cycle $Z(f)$.
\end{corollary}

\begin{proof}
By Theorem~\ref{KMTheorem} and Proposition~\ref{proppair}, we have 
\begin{align*}
\left\{\Lambda_{KM}(\eta),f\right\} &= a^+(0,0) \int_X \eta  \wedge e_q + 
\sum_{h\in L^\#/L} \sum_{n < 0}  a^+(h,n) \int_{Z(h,n)} \eta \\ 
 &=a^+(0,0) \big(\eta, e_q\big)_X +  \int_{Z(f)} \eta.
\end{align*}
The claim now follows from Theorem~\ref{Main2}. The above calculation also holds for $\eta = e_q^{p-1}$, see \cite{KMCan}; hence the corollary is valid in this case as well. $\Lambda_{B,\psi}(f)$ is harmonic by Remark~\ref{harmonic} or alternatively by \cite{KMCan}, Theorem~4.1.
\end{proof}


\begin{proof}[Proof of Theorem~\ref{Main1}]

We would like to use the adjointness of $R_{-k}$ and $-L_{2-k}$ to compute
\[
( \Theta(z,\varphi_{KM}), \xi_k(f))_{2-k,L} = ( \Theta(z,\varphi_{KM}), R_{-k}(v^k\bar{f}))_{2-k,L}.
\]
However, due to the rapid growth of $f$, the scalar product 
$( L_{2-k} \Theta(z,\varphi_{KM}), v^k\bar{f})_{-k ,L}$ does not converge. We have to proceed more carefully. 

Let $\mathcal{F}_t$ be  the truncated fundamental domain (see (\ref{Fdomain})) for the action of $\Sl_2(\Z)$ on $\h$. For $g \in A_{2-k,L}$ and $h \in A_{-k,L^-}$, we then have
(see \cite{Br1}, Lemma~4.2, correcting a sign error)
\begin{equation}\label{adjoint}
\int\limits_{\mathcal{F}_t} \langle g,R_{-k} h \rangle v^{-k} \,du \,dv =
 - \int\limits_{\mathcal{F}_t} \langle L_{2-k} g, h \rangle v^{2-k} \,du \,dv
\; + \; \int\limits_{-1/2}^{1/2} \langle g(u+it), h(u+it) \rangle v^{-k} \,du.
\end{equation}
Note that we do not require any regularity for $g$ or $h$ at the cusp.
We apply (\ref{adjoint}) for $g = \Theta(\tau,z,\varphi_{KM})$ and $h = v^k \bar{f}$ and obtain

\begin{lemma}\label{L1}
\begin{align}
( \Theta(z,\varphi_{KM}), R_{-k}(v^k\bar{f}))_{2-k,L} &= \lim_{t \to \infty} 
\int_{\mathcal{F}_t} \langle   \Theta(\tau,z,\varphi_{KM}), R_{-k} v^k\bar{f} \rangle v^{-k} \,du \,dv \label{A0} \\ 
&= \lim_{t \to \infty} 
\biggl( - \int_{\mathcal{F}_t}  \langle L_{2-k} \Theta(\tau,z,\varphi_{KM}), \bar{f} \rangle \frac{du\,dv}{v^2}  \label{A1}\\
 & \qquad \quad +  \int_{-1/2}^{1/2} \langle \Theta(u+it ,z,\varphi_{KM}), \overline{f(u+it)} \rangle \,du  \biggr) \label{A2}.
\end{align}
\end{lemma}

Next, we would like to compute the limit in (\ref{A1}) and (\ref{A2}) separately. We will see that this is possible for $z \notin Z(f)$, while for $z \in Z(f)$ the limits of both terms do not exist (in (\ref{A1}) and (\ref{A2}), the singularities cancel each other out). We first consider the term (\ref{A2}). The integration picks out the $0$-th Fourier coefficient of 
$\langle \Theta(u+it ,z,\varphi_{KM}), \overline{f(u+it)} \rangle$.

\begin{lemma}\label{L2}
Outside $Z(f)$, we have 
\begin{equation*}
\lim_{t \to \infty} \int_{-1/2}^{1/2} \langle \Theta(u+it ,z,\varphi_{KM}), \overline{f(u+it)} \rangle \,du =  a^+(0,0) e_q,
\end{equation*}
while the limit does not exist for $z \in Z(f)$. In particular, the limit defines a smooth differential form on $D-Z(f)$, which extends smoothly to $D$.
\end{lemma}

\begin{proof}

Since $f^-$ is rapidly decreasing in $v$, we see
\[
\lim_{t \to \infty} \int_{-1/2}^{1/2} \langle \Theta(u+it ,z,\varphi_{KM}), \overline{f^-(u+it)} \rangle \,du =0.
\]
We write $ \varphi_{KM}(x,z) = P_{KM}(x,z) \varphi_0(x,z)$, where $P_{{KM}}$ is a polynomial in 
$[\mathcal{P}(V(\R)) \otimes \mathcal{A}^{q}(D)]^G$ only depending on $x_{z^\perp}$. 
For $f^+$, one then has 
\begin{align*}
 \int_{-1/2}^{1/2} \langle \Theta(u+it & ,z,\varphi_{KM}), \overline{f^+(u+it)} \rangle \,du \\ 
&=
  \varphi_{KM} (0)  a^+(0,0) \; + \; \sum_{\lambda \in L^{\#}-0 } P_{{KM}}(\sqrt{t} \la,z)  a^+(\la,-q(\la)) e^{4 \pi q(\la_z) t}.
 \end{align*}
For $z \in Z(f)$, hence $\la_z =0$, the sum clearly diverges as $t \to \infty$, otherwise, the sum is rapidly decreasing in $t$ by arguments similar to the proof of Proposition~\ref{sing}. The lemma now follows from $\varphi_{KM}(0) = e_q$, see (\ref{Euler}).
\end{proof}


\begin{lemma}\label{L3} The quantity in (\ref{A1}) is given by
\begin{align*}
\lim_{t \to \infty} 
  \int_{\mathcal{F}_t}  \langle L_{2-k} \Theta(\tau,z,\varphi_{KM}), \bar{f} \rangle \,d\mu &= \int^{reg}_{\G' \back \h}  \langle L_{2-k} \Theta(\tau,z,\varphi_{KM}), \bar{f} \rangle\,d\mu\\
&= \int^{reg}_{\G' \back \h}  \langle \Theta(\tau,z, d \psi), \bar{f} \rangle \,d\mu.
\end{align*}
This defines a smooth form on $D-Z(f)$ which extends smoothly to all of $D$. 
\end{lemma}

\begin{proof}
The second equality follows from the fundamental fact
\[
 L_{2-k} \Theta(\tau,z,\varphi_{KM}) = \Theta(\tau,z, d \psi),
\]
which is Theorem~\ref{KMpsi}. From the Lemmas~\ref{L1} and \ref{L2} it follows that the limit exists, defining a smooth form outside $Z(f)$. As (\ref{A0}) defines a smooth form on $D$, the limit extends to all of $D$. As the limit exists, it is per definitionem the regularized integral. (This corresponds to $L \varphi_{KM}(0) = d \psi(0) = 0$). 
\end{proof}

Combining Lemmas~\ref{L1}, \ref{L2}, \ref{L3} we obtain outside $Z(f)$:
\begin{align}
( \Theta(z,\varphi_{KM}), \xi_k(f))_{2-k,L} = -  
\int^{reg}_{\G' \back \h}  \langle  \Theta(\tau,z, d \psi), \bar{f} \rangle \,d\mu \;  + \; a^+(0,0) e_q.
\end{align}
The identity (\ref{Main1a}) now follows from

\begin{lemma}\label{L4}
\[
\int^{reg}_{\G' \back \h}  \langle  \Theta(\tau,z, d \psi), \bar{f} \rangle \,d\mu = d 
\int^{reg}_{\G' \back \h}  \langle  \Theta(\tau,z,  \psi), \bar{f} \rangle \,d\mu = d \Phi(z,f,\psi).
\]
\end{lemma}

\begin{proof}
We recall that because of $P_\psi(0,z)=0$, we have
\[
\int^{reg}_{\G' \back \h}  \langle  \Theta(\tau,z,  \psi), \bar{f} \rangle \,d\mu =\lim_{t\to \infty}\int_{\mathcal{F}_t}  \langle  \Theta(\tau,z,\psi), \bar{f} \rangle \,d\mu .
\]
Since ${\mathcal{F}_t}$ is compact, we only have to 
justify the interchange of the differentiation $d$ with the limit 
$t\to\infty$. 
Arguing as in the proof of Proposition~\ref{sing} (using the notation from there) we find that it suffices to show:
\begin{align}\label{cc1}
d h(z,0) &=\sum_{\lambda\in L^\#}
a^+(\lambda,-q(\lambda))  \int_1^\infty d P_{\psi}(\sqrt{v}\la,z)\exp( 4\pi v q(\lambda_{z})) \frac{dv}{v}.
\end{align}
If $U\subset D$ is any relatively compact open subset and $\eps>0$, we split the sum over $\lambda\in L^\#$ into the finite sum over $\lambda\in S_f(U,\eps)$ and the sum over $\lambda\in L^\#-S_f(U,\eps)$.
An estimate for the latter sum as in Proposition~\ref{sing} shows that it converges uniformly absolutely on $U$, justifying the interchange of differentiation and the limit.
Thus it suffices to show that
\[
d \int_1^\infty P_{\psi}(\sqrt{v}\la,z)\exp( 4\pi v q(\lambda_{z})) 
\frac{dv}{v}
=\int_1^\infty d P_{\psi}(\sqrt{v}\la,z)\exp( 4\pi v q(\lambda_{z})) 
\frac{dv}{v}
\]
for $\lambda\in S_f(U,\eps)$ and $z\in U-Z(f)$. This follows from the exponential decay of the integrands as $v\to \infty$.
\end{proof}

For (\ref{Main1b}), we note that in this case, we have $\psi = -d^c \varphi_0$, see (\ref{dcformel}). Therefore, we just need to show
\begin{align}
\label{A6}
\int^{reg}_{\G' \back \h} d^c \langle \Theta(\tau,z,\varphi_0), \bar{f} \rangle \,d\mu 
&= d^c \int^{reg}_{\G' \back \h} \langle \Theta(\tau,z,\varphi_0), \bar{f} \rangle \,d\mu = d^c \Phi(z,f,\varphi_0). 
\end{align}
This goes through as above with one extra point: Since $\varphi_0(0)=1 \ne 0$, the constant term of $\langle \Theta(\tau,z,\varphi_0), \bar{f} \rangle$ involves also the term $ a^+(0,0)v$. One sees that the regularized integral is actually given by
\begin{equation}\label{Regul}
\int^{reg}_{\G' \back \h} \langle \Theta(\tau,z,\varphi_0), \bar{f} \rangle \,d\mu = a^+(0,0) C+ 
\lim_{t \to \infty} \int_{\mathcal{F}_t} \left(
\langle \Theta(\tau,z,\varphi_0), \bar{f} \rangle - a^+(0,0)v\right) \,d\mu, 
\end{equation}
where $C$ is a constant not depending on $z$ (namely the constant term in the Laurent expansion at $s=0$ of $\lim_{t \to \infty} \int_{\mathcal{F}_t} v^{1-s}d\mu$).
Applying $d^c$ to the right hand side, one can now interchange differentiation and the limit as above. One obtains \eqref{A6},
which implies (\ref{Main1b}).

This completes the proof of Theorem~\ref{Main1}.
\end{proof}


\section{The Borcherds lift as a current}

In this section, we consider the current which is induced by the Borcherds lift  $\Phi(f,\psi)$ defined in \eqref{genborch}.


Consider a top degree form $\phi \in A^{pq}(\G \back D)$. We then have
\[
\int_{\G \back D} \phi(z) = \left( \int_{\G \back G} \phi(g)dg \right) ({1}_{\mathfrak{p}}), 
\]
where $\phi$ on the right hand side is considered as an element in 
$[C^{\infty}(\G \back G) \otimes \bigwedge^{pq}\mathfrak{p}^{\ast}]^K$. (We will frequently use this identification without further comment.)
Here 
$1_{\mathfrak{p}}$ is a properly oriented basis vector for 
$\bigwedge^{pq}\mathfrak{p} \simeq \R$ of length one with respect to the Killing form. Moreover, $dg = dz \,dk$, where $dz$ is the measure on $D$ coming form the Killing form and $\vol(K,dk) =1$.

We now pick appropriate coordinates for $D$. 
We denote by $H=G_{v_1}$ the stabilizer in $G$ of the first basis vector $v_1$ of $V$. Then $H$ is the fixed point set of an involution $\tau$ on $G$. On the Lie algebra level we obtain a decomposition $\mathfrak{g} = \mathfrak{h} + \mathfrak{q}$ into $\pm1$ eigenspaces of $\tau$. Then $L = H \cap K$
 is a maximal compact subgroup of $H$. We write $\G_H = \G_{v_1}$.
There is a diffeomorphism (see \cite{KM2}, section 4, \cite{FJ}, section 2)
\begin{gather*}
\Psi: H \times _L(\mathfrak{p} \cap \mathfrak{q}) \longrightarrow  D=G/K, \\
(h,Y) \mapsto h \exp(Y) K.
\end{gather*}
Let $a_t =  \exp(tX_{1p+q})$ for $t\in \R$, and let $A = \{a_t ;\; t \in \R \} \simeq \R$ be the one-parameter subgroup of $G$ associated to the maximal abelian subspace of $\mathfrak{p} \cap \mathfrak{q}$ generated by $X_{1 p+q}$. We write $A_{\eps} = \{a_t; \;t \geq \eps\}$.   
We have a decomposition $G = HAK$ and, with a positive constant $C$ depending on the normalizations of the invariant measures, the integral formula (see \cite{FJ}, section 2)
\[
\int_G \phi(g) dg = C \int_K \int_{A_0} \int_H \phi(ha_t k) |\sinh(t)|^{q-1} \cosh(t)^{p-1} dh \, dt \, dk.
\]
(If $q=1$, one has to replace $A_0$ by $A$, corresponding to the fact that in this case $D_{v_1}$ splits $D$ into two pieces). 


\begin{proposition} 
The $(q-1)$-form $\Phi(f,\psi)$ is locally integrable, i.e.,
\[
\int_X \eta \wedge \Phi(f,\psi) < \infty
\]
for $\eta \in 
A_c^{(p-1)q+1}(X)$.
\end{proposition}

\begin{proof} 
For $q=1$ the statement is clear since in that case  $\Phi(f,\psi)$ is locally bounded, see Proposition~\ref{sing} and Example~\ref{hyper}. For $q \ne 1$, we need to show
\[
\int_{\G \back G}  \eta(g) \wedge \Phi(f,\psi)(g)  \,dg < \infty.
\]
We can replace $\Phi(f,\psi)(g)$ by the part which gives rise to its singularities 
\[
\sum_h \sum_{n \in \Q_{<0}} a^+(h,n) \sum_{\lambda \in L_{h,-n}} \tilde{\psi}(g^{-1}\lambda)
\]
with $\tilde{\psi}(\lambda) := \tilde{\psi}(\la, z_0,0)$ given by setting $s=0$ in (\ref{sing1}). 
%
Now
\[
\sum_{\lambda \in L_{h,-n}} \tilde{\psi}(g^{-1}\lambda) = \sum_{  \lambda \in \G \back L_{h,-n}} \sum_{\g \in \G_{\la} \back \G } 
 \tilde{\psi}(\g^{-1} g^{-1}\lambda),
\]
and therefore by unfolding
\begin{align*}
& \int_{\G \back G}  \sum_h \sum_{n \in \Q_{<0}} a^+(h,n) \sum_{\lambda \in L_{h,-n}} \eta(g) \wedge \tilde{\psi}(g^{-1}\lambda)  \,dg\\
&= 
\sum_h \sum_{n \in \Q_{<0}} a^+(h,n) 
\sum_{  \lambda \in \G \back L_{h,-n}}
\int_{\G_{\lambda} \back G}  \eta(g) \wedge \tilde{\psi}(g^{-1}\la) \,dg.
\end{align*}
This is a finite sum (since $a^+(h,n) =0$ for $n \ll 0$), so that we only have to show the existence of the integral for a given $\lambda$ of positive length. We can assume $\lambda = \sqrt{m}v_1$ for some $m >0$. We have 
\begin{multline}\label{INT}
\int_{\G_{H} \back G}  \eta(g) \wedge \tilde{\psi}(g^{-1}\sqrt{m}v_1) \,dg\\ 
= C' \int_{\G_{H} \back H} \int_0^{\infty}  \eta(ha_t) \wedge \tilde{\psi}( a_t^{-1} h^{-1} \sqrt{m}v_1) \sinh(t)^{q-1} \cosh(t)^{p-1} \,dt \,dh
\end{multline}
with a positive constant $C'$. But now $h^{-1} v_1 =v_1$ and 
\[
a_t^{-1} \sqrt{m} v_1 = \cosh(t) \sqrt{m} v_1 - \sinh(t) \sqrt{m} v_{p+q}.
\]
Hence
\begin{subequations}
\label{A-action}
\begin{align}\label{A-action1}
( a_t^{-1}\sqrt{m} v_1)_{z_0} &= -\sinh(t) \sqrt{m} v_{p+q}, \\
\label{A-action2}
( a_t^{-1} \sqrt{m} v_1)_{z_0^{\perp}} &= \cosh(t) \sqrt{m} v_1,
\end{align}
\end{subequations}
so that, see (\ref{sing1}),
\begin{multline}\label{psitilde}
\tilde{\psi}(a_t^{-1} \sqrt{m} v_1) = \frac{-1}{\sqrt{2\pi}}\sum_{\underline{\alpha}} R_{\underline{\alpha}}
(v_{p+q})
\sum_{\ell=0}^{q-1} Q_{\underline{\alpha},\ell}
\left(\sqrt{m}\cosh(t)v_1\right)\\ 
\times \left( 2\pi m\sinh^2(t) \right)^{-\ell/2} \G\left(\tfrac{\ell+1}{2}, 2\pi m \sinh^2(t)\right).
\end{multline}
Therefore the integrand in (\ref{INT}) is bounded as $t \to 0$. On the other hand, $\tilde{\psi}( a_t^{-1} h^{-1}v_1)$ is exponentially decreasing in $e^t$ (uniformly in $h$). Since $\eta$ has compact support and ${\G_{H} \back H}$ has finite volume, we conclude that (\ref{INT}) converges. Hence the above unfolding is valid and the proposition is proved.
\end{proof}

Recall that a locally integrable $\ell$-form $\omega$ on $X$ defines a current $[\omega]$ on $X$, i.e., a linear functional on $ \mathcal{A}^{pq-\ell}_c(X)$, via
\[
[\omega](\eta) := \int_{X} \eta \wedge \omega
\]
for $\eta \in  \mathcal{A}^{pq-\ell}_c(X)$. For a current $T$, we define its exterior derivative by
\[
(dT)(\eta) = (-1)^{\deg(\eta)+1} T(d\eta).
\]

\begin{theorem}
Let $\delta_{Z(f)}$ be the delta-current for the special cycle $Z(f)$, and  for a locally integrable differential form $\omega$, we denote the associated current by $[\omega]$. Then
\begin{align}\label{Haupt}
d[\Phi(f,\psi)] + \delta_{Z(f)} &= [\Lambda_{B,\psi}(f)]. 
\end{align}
\end{theorem}

\begin{proof}
Let $\eta \in \mathcal{A}_c^{(p-1)q}(X)$. Then
\begin{align*}
d[\Phi(f,\psi)](\eta) &= (-1)^{(p-1)q+1} \int_X d \eta \wedge \Phi(f,\psi)  \\ &= (-1)^{(p-1)q+1}\int_X d \left( \eta \wedge \Phi(f,\psi) \right) + \int_X \eta \wedge d \Phi(f,\psi) . 
\end{align*}
To obtain the theorem we therefore only need to show
\[
(-1)^{(p-1)q+1} \int_X d \left( \eta \wedge \Phi(f,\psi) \right) =  \int_{Z(f)} \eta.
\]
By Stokes' theorem, we can modify the differential form $\Phi(f,\psi)$ by a smooth form without changing the integral. Hence we can replace it by its `singular' part as in the  previous proposition. By unfolding we see
\begin{align*}
\int_X d \left( \eta \wedge \Phi(f,\psi)  \right) = 
\sum_h \sum_{n \in \Q_{<0}} a^+(h,n) 
\sum_{ \lambda \in \G \back L_{h,-n}} \int_{
\G_{\lambda} \back D} d \left( \eta \wedge \tilde{\psi}(\la,z,0) \right). 
\end{align*}
Again, this is a finite sum, so it suffices to consider the integral on the right hand side for $\la = \sqrt{m}v_1$.


For $\eps >0$,
\[
U_{\eps} :=  D - \Psi\left( H \times \left\{ tX_{1 p+q};\;t\geq \eps\right\} \right)
\]
defines an open neighborhood of the cycle $D_{v_1}$. (For $q=1$, replace the condition $t\geq \eps$ by $|t|\geq \eps$.) Then 
by Stokes' theorem we obtain
\begin{equation}\label{Stokes}
\int_{\G_{H} \back D} d \left( \eta \wedge \tilde{\psi}(\sqrt{m}v_1,z,0) \right)= \lim_{\eps \to 0} \int_{\partial\left( \G_H \back (D - U_{\eps}) 
\right)}  \eta \wedge \tilde{\psi}(\sqrt{m}v_1,z,0).
\end{equation}
By the analogue for (\ref{INT}) (which follows from the considerations in \cite{FJ}, section 2), we see that (\ref{Stokes}) is equal to
\begin{equation}\label{Stokes1}
C \lim_{\eps \to 0} \int_{\G_H \back H} 
\eta(h a_{\eps}) \wedge 
\tilde{\psi}(a_{-\eps} \sqrt{m}v_1)  \sinh(\eps)^{q-1} \cosh(\eps)^{p-1} \,dh \left(1_{\mathfrak{p} / \R X_{1p+q}}   \right).
\end{equation}
for some universal constant $C\neq 0$.

We consider (\ref{A-action}) and (\ref{psitilde}) (with $t = \eps$). For (\ref{Stokes1}) only the terms  with $\ell = q-1$ can contribute. But
\[
Q_{\underline{\alpha},q-1}(\sqrt{m}\cosh(\eps) v_1,z_0) = \begin{cases} 
\left( \sqrt{2m}\cosh(\eps) \right)^{q-1}, &\text{  $\underline{\alpha}=(1,1,\dots,1)$,} \\
0, &\text{otherwise.}
\end{cases}
\]
We obtain
\begin{equation*}
C \lim_{\eps \to 0} \int_{\G_H \back H}  \eta(h a_{\eps}) \wedge \tilde{\psi}(a_{-\eps} \sqrt{m}v_1) \sinh(\eps)^{q-1} \cosh(\eps)^{p-1} \,dh = C' \int_{\G_H \back H} \eta(h) \,dh
\end{equation*}
with a constant $C$ and $\tilde{\psi}(\lambda) = \tilde{\psi}(\la, z_0,0)$ as before. Therefore, the theorem holds with $\delta_{Z(f)}$ replaced by $C''\delta_{Z(f)}$  with a certain constant $C''$ independent of $\eta$. We conclude that the constant is equal to $1$
by noting that this is the case for $\eta$ closed: In (\ref{Haupt}),
$d[\Phi(f,\psi)](\eta) =0$, while $[\Lambda_{B,\psi}(f)](\eta) = (\eta,\Lambda_{B,\psi}(f))_X = \int_{Z(f)} \eta$ by Corollary~\ref{c99}.
\end{proof}

The previous theorem can be reformulated saying that the pair $\left(Z(f), \Phi(z,f,\psi) \right)$ defines a class in the group of differential characters introduced by Cheeger and Simons, see \cite{C,CS}. Note that we can consider 
the group of differential characters as the analogue for real manifolds of the arithmetic Chow group, see \cite{GS}.

\begin{theorem}
In the Hermitian case of $q=2$, then the `classical' Borcherds lift $\Phi(z,f)$ defines a Green's current for the divisor $Z(f)$, i.e.,
\[
dd^c[\Phi(z,f)] + \delta_{Z(f)} = 
 [\La_B(z,f)].
\]
\end{theorem}

\begin{proof}
Using $-d^c \Phi(z,f) =  \Phi(z,f,\psi)$, one easily obtains
\[
  (dd^c \eta) \Phi(z,f) = \eta \wedge dd^c \Phi(z,f) + d\left ((d^c \eta) \Phi(z,f)  +\eta  \wedge \Phi(z,f,\psi) \right).
\]
The theorem now follows from the previous theorem and the logarithmic growth of $\Phi(z,f)$ along $Z(f)$.
\end{proof}
 
This result also follows from the splitting \eqref{splitequ} of $\Phi(z,f)$, obtained in \cite{Br1} via the Fourier expansion, and the Poincar\'e-Lelong formula.

Observe that in the generic case the space $X$ will be non-compact. Thus, to study intersection theory, one has to consider a smooth compactification $\widetilde{X}$ of $X$. In particular, in the arithmetic intersection theory of special divisors at the Archimedian place, one has to study $\Phi(z,f)$ as a current on $\widetilde{X}$. This requires a detailed analysis of the growth of $\Phi(z,f)$ at the boundary, which turns out to be of $\log$- and $\log$-$\log$-type (for the case of Hilbert modular surfaces see \cite{BBK}). One has to work with the extension of arithmetic intersection theory provided in \cite{BKK}.

Oda and Tsuzuki \cite{OT} constructed Green's functions for the special divisors by means of Poincar\'e series 
involving the secondary spherical function.
By Theorem~4.7 of \cite{BK}, $\Phi(z,f)$ can be written as a linear combination of the Green's functions of 
 \cite{OT} in an explicit way.

Let us also compare $\Phi(z,f)$ with the Green's functions constructed by Kudla \cite{K, Ku2}.
He introduces for $x \in V(\R), x\ne 0$ and $z\in D$, the function
\begin{equation*}
\xi^0(x,z) = - \Ei(2\pi (x_z,x_z)),
\end{equation*}
which for $(x,x) >0$ turns out to be a Green's function for the cycle $D_x$ on the Hermitian domain $D$ and is smooth otherwise. Here, for $ w \in \C$, $\Ei(w) = \int_{-\infty}^w \tfrac{e^t}{t}dt$ is the exponential integral. While $\xi^0$ is \emph{not} a Schwartz function on $V(\R)$, we can still define $\xi(x,\tau,z) = \xi^0(\sqrt{v}x,z) e^{-\pi \tau (x,x)}$ for $\tau \in \h$, and one then has, see \cite{K,Ku2},
\[
dd^c \xi(x,\tau,z) = \varphi_{KM}(x,\tau,z)
\]
for $x \ne 0$.
On the other hand, we can apply the lowering operator $L_{\kappa}$ on $\h$ to $\xi$, and one easily checks
\[
L_{\kappa} \xi(x,\tau,z) = -\varphi_0(x,\tau,z).
\]
This provides a different proof of Theorem~\ref{ddcGleichung}, i.e.,
 $\phi := L_{\kappa}\varphi_{KM} = - dd^c\varphi_0$, 
and in summary, we obtain the following commutative diagram:
\[
\xymatrix{
 \xi(x,\tau,z) \;\ar @{|->}[r]^{L_{\kappa}} \ar @{|->}[d]^{dd^c} &  -\varphi_0(x,\tau,z)  \ar @{|->}[d]^{dd^c}\\
\varphi_{KM}(x,\tau,z)  \;\ar @{|->}[r]^{L_{\kappa}} &  \phi(x,\tau,z)}.
\]
One can then use the ideas of the proof of Theorem~\ref{Main1} together with
(\ref{Regul}) to account for the problem that $\xi$ is not defined for $x =0$, to express $\Phi(z,f)$ in terms of $\xi$. 

Finally, we would like to mention that one can use the ideas of this paper to give a somewhat different proof of the results in \cite{K} and  \cite{BK}  relating the geometric degrees (and the Archimedean contribution to the arithmetic degrees) of the special cycles to the coefficients of certain Eisenstein series and its derivatives (see also \cite{Ku2, KRY}). 

\end{document}